\documentstyle[12pt,amstex,amscd,graphicx]{amsart}

\begin{document}

\newtheorem{theorem}{Theorem}[section]
\newtheorem{prop}[theorem]{Proposition}
\newtheorem{lemma}[theorem]{Lemma}
\newtheorem{cor}[theorem]{Corollary}
\newtheorem{defn}[theorem]{Definition}

\title{Cannon-Thurston Maps for
 Pared Manifolds of Bounded Geometry}

\author{Mahan Mj }

\address{ RKM Vivekananda University \\ 
Belur Math, WB-711202, India \\
Email: mahan.mj@@gmail.com}
\thanks{Research partly supported
by a UGC Major Research Project  grant}
\date{}

\maketitle

\begin{center}
{\it To S.V. on his birthday} 
\end{center}

\begin{abstract}
Let $N^h \in H(M,P)$ be a hyperbolic structure of bounded geometry
on a  pared manifold
such that each component of $\partial_0 M = \partial M  -  P$ is
incompressible. We show that the
limit set of $N^h$ is locally connected by constructing a natural
Cannon-Thurston map. This provides a unified treatment, an alternate
proof and a  generalization of results due to Cannon and
Thurston, Minsky, Bowditch, Klarreich and the author.

\smallskip

\begin{center}


\end{center}

\end{abstract}

\tableofcontents

\section{Introduction}

This is one
 in a series of papers (in between \cite{mitra-trees} and
\cite{brahma-ibdd}) leading to
the existence of  Cannon-Thurston maps for, and local
connectivity of limit sets of, finitely generated  Kleinian groups. The project is completed in \cite{mahan-amalgeo} 
and \cite{mahan-split}. The main aim of this
 paper is to develop a
{\it reduction technique}. Given the
existence of Cannon-Thurston
maps for closed surface groups of bounded
geometry (cf. \cite{mitra-trees}), this paper develops techniques
to generalize this first to punctured surfaces
and then to bounded
geometry
pared manifolds with incompressible boundary. 
While this is not a reduction theorem per se, the {\it techniques} of this paper shall be
used mutatis mutandis in \cite{brahma-ibdd}
and \cite{mahan-split} to obtain generalizations
of Cannon-Thurston
theorems for surface groups to theorems for
pared manifolds with incompressible boundary.
The main focus of \cite{brahma-ibdd}, \cite{mahan-amalgeo}
and \cite{mahan-split} will be to describe 
geometries of closed surface groups
for which one can prove the existence of 
Cannon-Thurston maps.
The main theorem of this paper is:

\smallskip

{\bf Theorem \ref{main2}}:
Suppose that $N^h \in H(M,P)$ is a hyperbolic structure of {\em
  bounded geometry} 
on a pared manifold $(M,P)$ with incompressible boundary 
$\partial_0 {M} = (\partial M - P)$ (i.e. each component of
$\partial_0 M$ is incompressible). Let
$M_{gf}$ denote a geometrically finite hyperbolic structure adapted
to $(M,P)$. Then the map  $i: \widetilde{M_{gf}}
\rightarrow \widetilde{N^h}$ extends continuously to the boundary
$\hat{i}: \widehat{M_{gf}}
\rightarrow \widehat{N^h}$. If $\Lambda$ denotes the limit set of
$\widetilde{M}$, then $\Lambda$ is locally connected.

\medskip

Examples to which Theorem \ref{main2} applies include: \\

\begin{enumerate}
\item The cover corresponding to the fiber subgroup of a closed
  hyperbolic 3-manifold fibering over the circle (Cannon and Thurston
  \cite{CT}, now published as \cite{CTpub}). \\
\item Hyperbolic 3 manifolds of bounded geometry, which correspond to
  simply or doubly degenerate Kleinian groups isomorphic to closed
  surface groups (Minsky \cite{minsky-jams}). \\
\item Hyperbolic 3 manifolds of bounded geometry without parabolics
 and freely
  indecomposable fundamental group (Mitra - Section 4.3 of
 \cite{mitra-trees}, and Klarreich \cite{klarreich}). \\
\item Hyperbolic 3 manifolds of bounded geometry, arising from
  simply or doubly degenerate Kleinian groups corresponding to punctured
  surface groups (Bowditch \cite{bowditch-ct} \cite{bowditch-stacks}). \\
\end{enumerate}

The main issue  addressed in this paper has been raised in various
forms by Minsky \cite{minsky-cdm}, McMullen \cite{ctm-locconn} and the
author \cite{mitra-thesis}, \cite{bestvinahp}.

Let $M$ be a closed hyperbolic 3-manifold fibering over the circle with 
fiber $F$. Let $\widetilde F$ and $\widetilde M$ denote the universal
covers of $F$ and $M$ respectively. Then $\widetilde F$ and $\widetilde M$
are quasi-isometric to ${\Bbb{H}}^2$ and ${\Bbb{H}}^3$ respectively. Now let
${{\Bbb{D}}^2}={\Bbb{H}}^2\cup{\Bbb{S}}^1_\infty$ and 
${{\Bbb{D}}^3}={\Bbb{H}}^3\cup{\Bbb{S}}^2_\infty$
denote the standard compactifications. In \cite{CT} (now published as \cite{CTpub}) Cannon and Thurston
show that the usual inclusion of $\widetilde F$ into $\widetilde M$
extends to a continuous map from ${\Bbb{D}}^2$ to ${\Bbb{D}}^3$.

\medskip
This was generalized by Minsky \cite{minsky-jams},  and Klarreich 
\cite{klarreich} and
independently (and using different techniques) by the author
(Section 4 of \cite{mitra-trees})  to 
prove that if $M$ is a geometrically tame hyperbolic 3-manifold with
freely indecomposable fundamental group and injectivity 
radius bounded below, and
 if $M_{cc}$ denotes the (topological) compact
 core of 
$M$ then   the inclusion of ${\widetilde{M}}_{cc}$ into $\widetilde M$ extends to a 
continuous map from $\widehat M_{cc}$ to $\widehat M$  where  ${\widetilde M}_{cc}$ and 
 $\widetilde M$ denote the universal covers of $M_{cc}$ and $M$
respectively and $\widehat M_{cc}$
and $\widehat M$ denote the (Gromov) compactifications of $\widetilde M_{cc}$
and $\widetilde M$ respectively. However, all these results left unanswered 
the case of manifolds with parabolics. 

In \cite{ctm-locconn} McMullen proved 
that the corresponding result holds for punctured torus groups, using the model
 manifold built by Minsky for these groups in \cite{minsky-torus}. 
In \cite{brahma-ibdd}, we shall unify the framework of this paper with
a certain notion of {\it i-bounded geometry} to give a simultaneous
generalization of the results of this paper and those of McMullen
\cite{ctm-locconn}. 

The present paper was born of an attempt to find a new proof (along
the lines of \cite{mitra-trees}) of a
result of Bowditch. 
In 
\cite{bowditch-ct} \cite{bowditch-stacks},  Bowditch
proved the existence of Cannon-Thurston maps for punctured surface groups
 of bounded geometry using some of the ideas from the first paper in
 the present series by the author \cite{mitra-trees}. 
The main result of this paper simultaneously generalizes
Bowditch's results  \cite{bowditch-ct} \cite{bowditch-stacks} and
\cite{mitra-trees} and hence  includes the results of \cite{CT},
\cite{minsky-jams}, 
\cite{klarreich}, 
\cite{bowditch-ct}, \cite{bowditch-stacks} and
\cite{mitra-trees}. When $(M,P)$ is the 
pair $S \times I, \delta S \times I$, for $S$ a surface group, we get
Bowditch's 
result \cite{bowditch-ct} \cite{bowditch-stacks}. However, even in this case,
our proof is different and circumvents the use of locally infinite
(Farey) graphs 
as in the proofs of
\cite{bowditch-ct}  
\cite{bowditch-stacks}. 

We are grateful to Brian Bowditch for pointing out that local
connectivity (the second statement in Theorem \ref{main2} above)
also follows by combining Bowditch's result \cite{bowditch-ct} with a
result of Anderson-Maskit \cite{and-mask}. See Section 10 of
\cite{bowditch-ct} for details. The first part of Theorem \ref{main2}
answers a question attributed to Thurston by Bowditch in the same
paper (Section 10 of \cite{bowditch-ct}. See also \cite{abikoff}).

In a sense, Theorem \ref{main2} above is a
direct generalization of the following which we proved in
\cite{mitra-trees}

\medskip

{\bf Theorem:} (Theorem 4.7 of \cite{mitra-trees} )
 Let $M$ be a compact 3-manifold with 
incompressible
boundary $\partial {M}$ without torus components. 
 Let $M_{hyp}$ be a hyperbolic structure without parabolics 
on $M$ such  
 that each end of
the manifold has bounded geometry. Also suppose that $M_{gf}$ is a
geometrically finite structure on $M$. Then  \\
$\bullet$    the inclusion of $\widetilde M_{gf}$ into $\widetilde
  M_{hyp}$ 
extends to a 
continuous map from $\widehat M_{gf}$ to $\widehat M_{hyp}$ \\
$\bullet$ the limit set of   $\widetilde
  M_{hyp}$ is locally connected \\

\smallskip

This problem is a part of a more general problem.
 A natural question seems to be the following:  \\

{\bf Question:}\cite{mitra-trees} \cite{mitra-thesis} \cite{bestvinahp} 
Suppose $H$ is a hyperbolic group acting freely and properly discontinuously 
by isometries on a hyperbolic metric space $X$.  Does 
the continuous proper map $ i$ : $\Gamma_H \rightarrow X$
extend to a continuous map $\hat i$ : 
$\widehat{\Gamma_H} \rightarrow \widehat{X}$ ? 

We can ask the same question for
relatively hyperbolic groups (in the sense of Gromov\cite{gromov-hypgps},
 Farb\cite{farb-relhyp}
and Bowditch \cite{bowditch-relhyp}). A convenient framework for
formulating this question is that of convergence groups acting on
compacta (Bowditch \cite{bowditch-cgnce} ). Hyperbolic and relatively
hyperbolic groups have been characterized in this setting by Bowditch
\cite{bowditch-jams} and Yaman \cite{yaman-relhyp} .

{\bf Question:} Suppose a relatively hyperbolic group $H$ acts on a
compact set $K$ such that the action is a convergence action and such
that each `cusp group' of $H$ leaves a point of $K$ fixed. Does there
exist an equivariant continuous map from ${\partial}H$ to $K$? (Here
${\partial}H$ denotes the boundary of the relatively group $H$, which
is well-defined as per \cite{bowditch-relhyp}. )

A closely related question was raised by McMullen \cite{ctm-locconn} 
in the context of Kleinian 
groups. However, his
question deals with the Floyd 
completion rather than the Gromov completion as we 
have done here. There 
is some difference between the two notions for manifolds with
cusps. Nevertheless, 
modulo Floyd's result \cite{floyd}, the two questions are equivalent in 
the context of Kleinian groups. This question has also gained
attention after the recent resolution of the Ending Lamination
Conjecture by Minsky \cite{minsky-elc1}
and Brock-Canary-Minsky
\cite{minsky-elc2}. In fact the question in the context of Kleinian
groups
has also been raised by Minsky in
\cite{minsky-cdm}. Both McMullen and Minsky ask if the limit sets of
Kleinian groups are locally connected given that they are connected. This will be completely answered in the affirmative in \cite{mahan-split}.

 Perhaps the most general context in which the   question of
 local connectivity  makes sense
 is that of convergence group actions on compacta \cite{bowditch-cgnce}.

{\bf Question:} Suppose a finitely generated group $H$ acts on a
compact connected metrizable perfect set $K$ such that the action is a
convergence action. Is $K$ locally connected? (i.e. does `admitting a
convergence action' promote a connected compactum to a Peano continuum?)

\subsection{Punctured Surface Groups of Bounded Geometry: A New Proof}

As a motivational example, we sketch the proof for Kleinian groups $G$ that
correspond to punctured surface groups, such that the 3-manifold $N^h = 
{{\Bbb{H}}^3}/G$ has bounded geometry away from the cusps, i.e. closed
geodesics have length uniformly bounded away from zero. These are
precisely the examples handled by Bowditch in \cite{bowditch-ct} and
\cite{bowditch-stacks} . We give a sketch in this particular case
because these are the simplest new non-trivial examples and also to
underscore the difference in our approach. 

We first excise
the cusps (if any) of $N^h$ leaving us a manifold that has one or two ends. Let
$N$ denote $N^h$ minus cusps. Then $N$ is quasi-isometric to the
universal curve over a Lipschitz path in Teichmuller space from which
cusps have been removed. (In fact, in \cite{mitra-trees} we had proven
that the path in question is a
Teichmuller 
geodesic. The proof there was for closed surfaces, but can be extended
painlessly to that of surfaces with punctures). This path is semi-infinite
or bi-infinite 
according as $N$ is one-ended or two-ended. Fix a reference finite
volume hyperbolic surface $S^h$. Let $S$ denote $S^h$ minus 
cusps. Then $\widetilde{S}$ is quasi-isometric to the Cayley graph of 
$\pi_1{(S)}$ which is (Gromov) hyperbolic. We fix a base surface in $N$
and identify it with $S$. Now look at
$\widetilde{S}\subset\widetilde{N}$. 
 Let $\lambda = [a,b]$ be a geodesic segment in $\widetilde{S}$.
We `flow' $\lambda$ out the end(s) of $\widetilde{N}$ to generate
a {\bf hyperbolic ladder-like set} $B_{\lambda}$ (thinking of
$\widetilde{S} $ as a horizontal sheet, See Figure 4)
that looks 
topologically
 like $[a,b] \times [0,\infty )$ or $[a,b]\times ({-}\infty , \infty)$ 
according as $N$ has one or two ends. This is exactly
a reproduction 
of our construction in \cite{mitra-trees} or \cite{mitra-ct} . 

A few details are in order. Regarding
$N$ as a universal curve (minus cusps) over a Lipschitz path in Teichmuller
 space, we can assume
that each fibre over the path is topologically $S$ with hyperbolic structure
varying over the path. The union of all points that correspond to $a$ (or $b$)
is a quasigeodesic ray (or bi-infinite quasigeodesic)
 in  $\widetilde{N}$. Join the pairs of points
that lie in the lift of a single copy of
$\widetilde{S}\subset\widetilde{N}$,  giving
a geodesic in each copy of $\widetilde{S}$.  This is
what is meant by `flowing'
$\lambda$ out the end(s) with the two quasigeodesic rays mentioned
above as `guides' (or 
 boundary lines), generating $B_{\lambda}$.
 The main technical theorems of \cite{mitra-ct} (Theorem 3.7) or
 \cite{mitra-trees} 
(Theorem 3.8) ensure that there is a retraction from  $\widetilde{N}$ to
$B_\lambda$ which does not increase distances much. From this it
 follows that $B_\lambda$ is 
quasi-isometrically embedded in $\widetilde{N}$. [A brief proof of
  this last
  assertion, viz. q.i. embeddedness follows from retract, was given
  in \cite{mitra-ct} but we had omitted it in \cite{mitra-trees}. We
  are grateful to Brian Bowditch for pointing out this gap in
  \cite{bowditch-ct}. See Gromov \cite{gromov-hypgps}, Section 7.3J, p.197
  for related results, and \cite{minsky-bddgeom} for a detailed proof
  of the same.]  Note that for
this we do not need 
$\widetilde{N}$ to be hyperbolic (in fact $\widetilde{N}$ is hyperbolic only
when $S^h$ has no cusps and $\widetilde{N^h} = \widetilde{N}$. This is what we had used
to prove the existence of Cannon-Thurston maps in \cite{mitra-ct} and \cite{mitra-trees}). 

Now if $\lambda$ lies outside a large ball about a fixed reference point $p$ in 
$\widetilde{S}$, then so does $B_\lambda$ in $\widetilde{N}$.  Since
$B_\lambda$ is q.i. embedded in $\widetilde{N}$, there exists an 
ambient 
$\widetilde{N}$-quasigeodesic $\mu$ lying in a bounded neighborhood of 
$B_\lambda$ joining the end-points of $\lambda$. If $S^h$ had no cusps, we immediately
conclude that for any geodesic segment $\lambda$ in $\widetilde{S^h}$ lying outside
large balls around $p$, there is a quasigeodesic in $\widetilde{N^h}$ joining its
endpoints and lying outside a large ball around $p$  in
$\widetilde{N^h}$. This gives us a continuous 
extension of the inclusion of  $\widetilde{S^h}$ into $\widetilde{N^h}$ to the boundary. At this
stage we have recaptured the results of Cannon and Thurston \cite{CT} and Minsky
\cite{minsky-jams}. 

However, when $S^h$ has cusps, $\widetilde{S^h}$ and $\widetilde{S}$ are different. So a 
little more work is necessary. Suppose as before that $\lambda_0$ is a geodesic in
$\widetilde{S^h}$ lying outside a large ball around $p$. For ease of exposition we assume
that the end-points of $\lambda_0$ lie outside cusps.
 Let $\lambda \subset \widetilde{S}$ be the geodesic in (the path-metric on)
$\widetilde{S}$ joining the same pair of points. Then off horodisks,
$\lambda_0$ and $\lambda$ track each other. Construct $B_\lambda$ as before, and let
$\mu$ be an ambient $\widetilde{N}$-quasigeodesic lying in a bounded neighborhood of 
$B_\lambda$ joining the end-points of $\lambda$. Then off horoballs in 
$\widetilde{N^h}$, $\mu$ lies outside a large ball around $p$. Let $\mu_0$ be the
hyperbolic geodesic joining the end points of $\mu$. Off horoballs, $\mu$ and
$\mu_0$ track each other. Hence, off horoballs, $\mu_0$ lies outside large balls
about $p$. The points at which $\mu_0$ intersects a particular horoball therefore lie 
 outside large balls 
about $p$. But then the hyperbolic segment joining them must do the same. This shows
that $\mu_0$ must itself lie outside large balls around $p$. As before we conclude
that there exists  a continuous
extension of the inclusion of  $\widetilde{S^h}$ into $\widetilde{N^h}$ to
the boundary. At this 
stage we have recaptured the result of Bowditch for punctured surfaces \cite{bowditch-ct}
\cite{bowditch-stacks}.

We note here that the main purpose of this paper is to continue with
our parallel approach to handling problems like the existence of a
Cannon-Thurston map, and local connectivity of limit sets, 
initiated in \cite{mitra-ct} and
\cite{mitra-trees}. These papers and the current one circumvent much
of the sophisticated  machinery specific to 2- and 3 manifolds
developed by Thurston {\it et al}. Thus, they provide a direct
approach to the existence of a Cannon-Thurston map, (as also local
connectivity), without having to deal with laminations and associated
geometries. Of course, this means that we lose out on the explicit
description that Cannon and Thurston \cite{CT} or Minsky
\cite{minsky-jams} or Bowditch \cite{bowditch-ct} provide in terms of
laminations. McMullen \cite{ctm-locconn} also uses the notion of
foliations in his proof to derive explicit information about locations
of closed geodesics.
However, the large-scale or coarse nature of our  approach makes it 
more general, and suitable
for application to cases where the surface machinery of Thurston is
not available, e.g. general hyperbolic or relatively hyperbolic spaces
and groups in the sense of Gromov, and more specifically, free groups
where the notion of laminations developed by Bestvina, Feighn and
Handel is substantially different \cite{BH-tt}
\cite{BFH-lam}. However, once the Cannon-Thurston 
map is in place, some notions of laminations can be resurrected, as in
\cite{mitra-endlam}. Thus, in principle, a long term aim of the
present approach would be to derive further parallels between 3
manifolds and surfaces on the one hand, and discrete groups and their
subgroups on the other. The latter context being vastly more general,
we have tried to avoid the techniques specific to 2 and 3 dimensions. 

The main technical difficulty in the first part of this
 paper (up to Section 4)  stems from the fact that a closed
geodesic on a component of the boundary of the pared manifold
$\partial_0 M = (\partial M - P)$ may be homotopic to a curve on a
 component  of $P$, i.e.
 it is an accidental parabolic in any hyperbolic structure on
 $(M,P)$. 
This
results in a somewhat trying case-by-case analysis in Section 4
of this paper. The resulting notation becomes a bit elaborate at times
and we pause at a few moments during the course of the paper to
summarize notation and keep it clear. 
There are essentially 3 cases we need to handle: \\
$\bullet$ {\bf $Z$-cusps:} These are relatively easy to handle as the
above sketch shows. \\
$\bullet$ {\bf $(Z + Z)$ cusps}, where no curve on any component of
$\partial_0 M$ is homotopic to a curve on the boundary torus
corresponding to the  $(Z + Z)$ cusp. \\
$\bullet$ {\bf $Z$ and  $(Z + Z)$ cusps}, where some curve(s) on
 component(s) of  
$\partial_0 M$ is(are)  homotopic to a particular  curve on the boundary torus
corresponding to the  $(Z + Z)$ cusp or some multiple of the core curve
of the $Z$ - cusp.

We shall have occasion to introduce a certain extra hypothesis of {\it
  p-incompressibility}, but 
we stick to the general case, i.e. we {\it do not} introduce the extra
  hypothesis of {\it
  p-incompressibility} (or absence of accidental parabolics on
  components of $\partial_0 M$) till Section 5.2, because the construction of
  $B_\lambda$ goes through even in the absence of this simplifying
  assumption. Further, once the proof is completed for the {\it
  p-incompressible} case, we use it (towards the end of Section 5) to
  prove the general case.

\subsection{Outline of the paper}

A brief outline of the paper follows. Section 2 deals with preliminaries on 
(Gromov) hyperbolic spaces. We also recall a result from McMullen 
\cite{ctm-locconn} which
says roughly that hyperbolic geodesics and ambient geodesics (geodesics in spaces
obtained by removing some horoballs) track each other off horoballs. 

In Section 3, we recall Thurston's definition of pared manifolds \cite{thurston-hypstr1}
\cite{thurston-hypstr3} and show that the universal cover of
a hyperbolic 3-manifold $M$, whose compact core is
a pared manifold, is quasi-isometric to a tree $T$ of hyperbolic metric spaces with possibly
exceptional vertex, once the cusps of $M$ have been removed. Further, if $\alpha$ denote the
root vertex, we demand that $\alpha$
be the possibly exceptional vertex. Let $X_\alpha$ denote the
corresponding vertex space.We shall choose a geometrically finite structure $M_{gf}$
on $M$ and identify $X_\alpha$ with ($\widetilde{M_{gf}}$ minus
$Z$-cusps). This  converts
$X_\alpha$ into a hyperbolic metric space, but  we  relax
the requirement that the embeddings of the edge spaces 
into $X_\alpha$ be qi-embeddings ({\it This is what makes $\alpha$ exceptional}). 
All other vertex and edge spaces are
 hyperbolic and all other inclusions of edge
spaces into vertex spaces are qi-embeddings.

In Section 4, we modify the construction of \cite{mitra-trees} to
construct a quasi-isometrically embedded {\bf hyperbolic ladder-like}
set $B_\lambda$ (See Figure 4)
corresponding to a geodesic $\lambda$. 
 We  show that $B_\lambda$ is 
qi-embedded in ( $ \widetilde{M} -  Z$-cusps). $Z$-cusps and
$(Z+Z)$-cusps are treated differently.The construction of $B_\lambda$
does not require {\it p-incompressibility} of the boundary components,
but only their incompressibility. 

It is in Section 5, that we first restrict the scope to pared manifolds with
{\it p-incompressible} boundary components. We show that
 if $\lambda$ lies outside
a large ball about a fixed reference point $p$ in $(X_\alpha \cup
cusps)$
 modulo
horoballs then so does $B_\lambda$ in
$\widetilde{M}$.

Finally, we use the tracking properties of ambient quasigeodesics vis a vis
hyperbolic geodesics and assemble the proof of the main theorem in the
case of {\it p-incompressible boundary}. We
deduce from this that 
the limit sets of the corresponding Kleinian groups are locally connected.
Once this case is proven, we use it to prove the result for pared
manifolds of incompressible boundary, thus relaxing the assumption on
p-incompressibility. 
The concluding Section 6 deals with examples to which our theorem applies, notably
Brock's example \cite{brock-itn}. We also indicate possible directions of
generalization. It would be worth bearing in mind that most of
the arguments of this paper are relevant to the considerably more
general framework of relatively hyperbolic groups {\it a la} Gromov
\cite{gromov-hypgps}, Farb \cite{farb-relhyp}, Bowditch
\cite{bowditch-relhyp}.

\smallskip

{\bf Apologia} In {\it From Beowulf to Virginia Woolf: An Astounding and Wholly Unauthorized History of English Literature}
\cite{rmmyers} Robert Manson Myers claims that the plays usually attributed to Shakespeare
are not in fact written by him but by another person of the same
name. Being fully aware
of our literary deficiencies we frankly admit that the papers \cite{mitra-ct}, \cite{mitra-trees}
and the present one are written by the same person under different names.

\smallskip

{\bf Acknowledgements} I am grateful to Brian Bowditch for
several helpful comments on a previous version of this paper.
I would also like to thank the referee for carefully
 reading the manuscript and suggesting several corrections and improvements.

\section{Preliminaries}

We start off with some preliminaries about hyperbolic metric
spaces  in the sense
of Gromov \cite{gromov-hypgps}. For details, see \cite{CDP}, \cite{GhH}. Let $(X,d)$
be a hyperbolic metric space. The 
{\bf Gromov boundary} of 
 $X$, denoted by $\partial{X}$,
is the collection of equivalence classes of geodesic rays $r:[0,\infty)
\rightarrow\Gamma$ with $r(0)=x_0$ for some fixed ${x_0}\in{X}$,
where rays $r_1$
and $r_2$ are equivalent if $sup\{ d(r_1(t),r_2(t))\}<\infty$.
Let $\widehat{X}$=$X\cup\partial{X}$ denote the natural
 compactification of $X$ topologized the usual way(cf.\cite{GhH} pg. 124).

The {\bf Gromov inner product}
 of elements $a$ and $b$ relative to $c$ is defined 
by 
\begin{center}
$(a,b)_c$=1/2$[d(a,c)+d(b,c)-d(a,b)]$.
\end{center}

\begin{defn} A subset $Z$ of $X$ is said to be 
{\bf $k$-quasiconvex}
 if any geodesic joining points of  $ Z$ lies in a $k$-neighborhood of $Z$.
A subset $Z$ is {\bf quasiconvex} if it is $k$-quasiconvex for some
$k$. 
\end{defn}

For  simply connected real hyperbolic
manifolds this is equivalent to saying that the convex hull of the set
$Z$ lies in a  bounded neighborhood of $Z$. We shall have occasion to
use this alternate characterization.

\begin{defn}
A map $f$ from one metric space $(Y,{d_Y})$ into another metric space 
$(Z,{d_Z})$ is said to be
 a {\bf $(K,\epsilon)$-quasi-isometric embedding} if
 
\begin{center}
${\frac{1}{K}}({d_Y}({y_1},{y_2}))-\epsilon\leq{d_Z}(f({y_1}),f({y_2}))\leq{K}{d_Y}({y_1},{y_2})+\epsilon$
\end{center}
If  $f$ is a quasi-isometric embedding, 
 and every point of $Z$ lies at a uniformly bounded distance
from some $f(y)$ then $f$ is said to be a {\bf quasi-isometry}.
A $(K,{\epsilon})$-quasi-isometric embedding that is a quasi-isometry
will be called a $(K,{\epsilon})$-quasi-isometry.

A {\bf $(K,\epsilon)$-quasigeodesic}
 is a $(K,\epsilon)$-quasi-isometric embedding
of
a closed interval in $\Bbb{R}$. A $(K,0)$-quasigeodesic will also be called
a $K$-quasigeodesic.
\end{defn}

Let $(X,{d_X})$ be a hyperbolic metric space and $Y$ be a subspace that is
hyperbolic with the inherited path metric $d_Y$.
By 
adjoining the Gromov boundaries $\partial{X}$ and $\partial{Y}$
 to $X$ and $Y$, one obtains their compactifications
$\widehat{X}$ and $\widehat{Y}$ respectively.

Let $ i :Y \rightarrow X$ denote inclusion.

\begin{defn}  Let $X$ and $Y$ be hyperbolic metric spaces and
$i : Y \rightarrow X$ be an embedding. 
 A {\bf Cannon-Thurston map} $\hat{i}$  from $\widehat{Y}$ to
 $\widehat{X}$ is a continuous extension of $i$. 
\end{defn}

The following  lemma (Lemma 2.1 of \cite{mitra-ct})
 says that a Cannon-Thurston map exists
if for all $M > 0$ and $y \in Y$, there exists $N > 0$ such that if $\lambda$
lies outside an $N$ ball around $y$ in $Y$ then
any geodesic in $X$ joining the end-points of $\lambda$ lies
outside the $M$ ball around $i(y)$ in $X$.
For convenience of use later on, we state this somewhat
differently.

\begin{lemma}
A Cannon-Thurston map from $\widehat{Y}$ to $\widehat{X}$
 exists if  the following condition is satisfied:\\
Given ${y_0}\in{Y}$, there exists a non-negative function  $M(N)$, such that 
 $M(N)\rightarrow\infty$ as $N\rightarrow\infty$ and for all geodesic segments
 $\lambda$  lying outside an $N$-ball
around ${y_0}\in{Y}$  any geodesic segment in $X$ joining
the end-points of $i(\lambda)$ lies outside the $M(N)$-ball around 
$i({y_0})\in{X}$.
\label{contlemma}
\end{lemma}

The above result can be interpreted as saying that a Cannon-Thurston map 
exists if the space of geodesic segments in $Y$ embeds properly in the
space of geodesic segments in $X$.

\smallskip

We shall also be requiring certain properties of hyperbolic spaces
minus horoballs. These were studied by Farb \cite{farb-relhyp} under
the garb of `electric geometry'. We
combine Farb's results with a 
version that is a (slight variant of) theorem due to McMullen
(Theorem 8.1 of \cite{ctm-locconn}).

\begin{defn} A path $\gamma : I \rightarrow Y$ to a path metric space $Y$ is an ambient
K-quasigeodesic if we have
\begin{center}
$L({\beta}) \leq K L(A) + K$
\end{center}
for any subsegment $\beta = \gamma |[a,b]$ and any path $A : [a,b] \rightarrow Y$ with the
same endpoints. 
\end{defn}

The following definitions are adapted from \cite{farb-relhyp}

\begin{defn} Let $M$ be a convex hyperbolic manifold.
Let $Y$ be the universal cover of  $M$ minus cusps
and $X = \widetilde{M}$.  $\gamma$ is said to be a $K$-quasigeodesic in $X$ 
{ \bf without backtracking } if \\ 
$\bullet$ $\gamma$ is  a
  $K$-quasigeodesic in $X$ \\
$\bullet$  $\gamma$
does not return to any
  horoball $\bf{H}$ after leaving it. \\

\smallskip

 $\gamma$ is said to be an ambient
$K$-quasigeodesic in $Y$
{ \bf without backtracking } if \\ 
$\bullet$ $\gamma$ is an ambient $K$-quasigeodesic in $Y$ \\
$\bullet$ $\gamma$ is obtained from a $K$-quasigeodesic without
backtracking in $X$ by
replacing each  maximal subsegment with end-points on a horosphere by
a quasigeodesic lying on the surface of the horosphere. 
\end{defn}

Note that in the above definition, we allow the behavior to be quite
arbitrary on horospheres (since Euclidean quasigeodesics may be quite
wild); however, we do not allow wild behavior off horoballs.

$B_R (Z)$ will denote the $R$-neighborhood of the set $Z$. \\
Let $\cal{H}$ be a locally finite collection of horoballs in a convex
subset $X$ of ${\Bbb{H}}^n$ 
(where the intersection of a horoball, which meets $\partial X$ in a point, 
 with $X$ is
called a horoball in $X$). 
The following theorem is due to McMullen \cite{ctm-locconn}.

\begin{theorem} \cite{ctm-locconn} 
Let $\gamma: I \rightarrow X \setminus \bigcup \cal{H} = Y$ be an ambient
$K$-quasigeodesic for a convex subset $X$ of ${\Bbb{H}}^n$ and let
$\mathcal{H}$  denote a
collection of horoballs.
Let $\eta$ be the hyperbolic geodesic with the same endpoints as
$\gamma$. Let $\cal{H}({\eta})$  
be the union of all the horoballs in $\cal{H}$ meeting $\eta$. Then
$\eta\cup\mathcal{H}{({\eta})}$ is (uniformly) quasiconvex and $\gamma
(I) \subset  
B_R (\eta \cup \cal{H} ({\eta}))$, where $R$ depends only on
$K$. 
\label{ctm}
\end{theorem}

The above theorem is similar in flavor to certain theorems
 about relative hyperbolicity
{\it a la} Gromov \cite{gromov-hypgps}, Farb \cite{farb-relhyp} and Bowditch
\cite{bowditch-relhyp}. We give below a related theorem that is
derived from Farb's
`Bounded Horosphere Penetration' property. 

Let $\gamma_1 = \overline{pq}$ be a hyperbolic $K$-quasigeodesic
without backtracking starting from a horoball $\bf{H_1}$ and ending
within (or on) a {\em different}
horoball ${\bf{H_2}}$. Let $\gamma = [a,b]$ be the hyperbolic geodesic
minimizing distance between $\bf{H_1}$ and $\bf{H_2}$. Following
\cite{farb-relhyp} we put the zero metric on the horoballs that
$\gamma$ penetrates. The resultant pseudo-metric is called the
electric metric. Let $\widehat{\gamma}$ and  $\widehat{\gamma_1}$
denote denote the paths represented by $\gamma$ and $\gamma_1$
respectively in this pseudometric. It is shown in \cite{farb-relhyp}
that $\gamma$, $\widehat{\gamma}$ and  $\widehat{\gamma_1}$ have
similar intersection patterns with horoballs, i.e. 
there exists $C_0$ such that \\
$\bullet 1$ If only one of $\gamma$ and $\widehat{\gamma_1}$
 penetrates a horoball
  $\bf{H}$, then it can do so for a distance $ \leq C_0$. \\
$\bullet 2$ If both $\widehat{\gamma_1}$ and $\gamma$ enter (or leave) a horoball
  $\bf{H}$ then their entry (or exit) points are at a distance of at
  most $C_0$ from each other. [Here by `entry' (resp. `exit') point of a
  quasigeodesic we mean
  a point at which the path switches from being in the complement of
  or `outside'
 (resp. in the interior of or `inside') a closed horoball to being
  inside (resp. outside) it].\\
The point to observe here is that quasigeodesics without backtracking
in our definition gives rise to quasigeodesics without backtracking in
 Farb's sense. Since this is true for arbitrary $\gamma_1$ we give
 below a slight strengthening of this fact. Further, by our
 construction of ambient quasigeodesics without backtracking, we might
 just as well consider ambient quasigeodesics without backtracking in
 place of quasigeodesics.

\begin{theorem}  \cite{farb-relhyp}
Given $C > 0$, there exists $C_0$ such that if \\
$\bullet 1$ either two quasigeodesics without backtracking 
$\gamma_1 , \gamma_2$ in $X$, OR\\
$\bullet 2$ two ambient
quasigeodesics without backtracking  $\gamma_1 , \gamma_2$ in $Y$, OR \\
$\bullet 3$ $\gamma_1$ - an ambient
quasigeodesic without backtracking in $Y$ and $\gamma_2$ -
a quasigeodesic without
backtracking in $X$,\\

\smallskip

\noindent start and end \\

\smallskip

\noindent $\bullet 1$ either  on (or within) the same horoball OR\\
$\bullet 2$  a distance $C$ from each other \\
then they have similar intersection patterns with horoballs (except
possibly the first and last ones), i.e.
there exists $C_0$ such that \\
$\bullet 1$ If only $\gamma_1$ penetrates or travels along the boundary of a horoball
  $\bf{H}$, then it can do so for a distance $ \leq C_0$. \\
$\bullet 2$ If both $\gamma_1$ and $\gamma_2$ enter (or leave) a horoball
  $\bf{H}$ then their entry (or exit) points are at a distance of at
  most $C_0$ from each other. 
\label{farb}
\end{theorem}

\section{Trees of Hyperbolic Metric Spaces and Pared Manifolds}

\subsection{Definitions}

We start with a notion closely related to one  introduced in
\cite{BF}. 
\begin{defn}
 A  {\bf tree (T) of hyperbolic metric spaces satisfying
the q(uasi) i(sometrically) embedded condition} is a metric space $(X,d)$
admitting a map $P : X \rightarrow T$ onto a simplicial tree $T$, such
that there exist $\delta {,} \epsilon$ and $K > 0$ satisfying the following: \\
\begin{enumerate}
\item  For all vertices $v\in{T}$, 
$X_v = P^{-1}(v) \subset X$ with the induced path metric $d_v$ is a 
$\delta$-hyperbolic metric space. Further, the
inclusions ${i_v}:{X_v}\rightarrow{X}$ 
are uniformly proper, i.e. for all $M > 0$, 
there exists $N > 0$ such that for all 
$v\in{T}$ and $x, y\in{X_v}$,
 $d({i_v}(x),{i_v}(y)) \leq M$ implies
${d_v}(x,y) \leq N$.
\item Let $e$ be an edge of $T$ with initial and final vertices $v_1$ and
$v_2$ respectively.
Let $X_e$ be the pre-image under $P$ of the mid-point of  $e$.  
Then $X_e$ with the induced path metric is $\delta$-hyperbolic.
\item There exist maps ${f_e}:{X_e}{\times}[0,1]\rightarrow{X}$, such that
$f_e{|}_{{X_e}{\times}(0,1)}$ is an isometry onto the pre-image of the
interior of $e$ equipped with the path metric.
\item ${f_e}|_{{X_e}{\times}\{{0}\}}$ and 
${f_e}|_{{X_e}{\times}\{{1}\}}$ are $(K,{\epsilon})$-quasi-isometric
embeddings into $X_{v_1}$ and $X_{v_2}$ respectively.
${f_e}|_{{X_e}{\times}\{{0}\}}$ and 
${f_e}|_{{X_e}{\times}\{{1}\}}$ will occasionally be referred to as
$f_{v_1}$ and $f_{v_2}$ respectively.
\end{enumerate}   
\end{defn}

$d_v$ and $d_e$ will denote path metrics on $X_v$ and $X_e$ respectively.
$i_v$, $i_e$ will denote inclusion of $X_v$, $X_e$ respectively into $X$.

We need a version of the above definition adapted to 3 manifolds with cusps.
For convenience of exposition, $T$ shall be assumed to be rooted, i.e.
equipped with a base vertex $\alpha$.

\begin{defn}
A  tree (T) of hyperbolic metric spaces with {\bf possibly
exceptional vertex} satisfying
the q(uasi) i(sometrically) embedded condition is a metric space $(X,d)$
admitting a map $P : X \rightarrow T$ onto a rooted simplicial tree
$T$ with root $\alpha$, such
that there exist $\delta{,} \epsilon$ and $K > 0$ satisfying the following: \\
\begin{enumerate}
\item  For all vertices $v\in{T}$, 
$X_v = P^{-1}(v) \subset X$ with the induced path metric $d_v$ is a 
$\delta$-hyperbolic metric space. Further, the
inclusions ${i_v}:{X_v}\rightarrow{X}$ 
are uniformly proper, i.e. for all $M > 0$, $v\in{T}$ and $x, y\in{X_v}$,
there exists $N > 0$ such that $d({i_v}(x),{i_v}(y)) \leq M$ implies
${d_v}(x,y) \leq N$.
\item Let $e$ be an edge of $T$ with initial and final vertices $v_1$ and
$v_2$ respectively.
Let $X_e$ be the pre-image under $P$ of the mid-point of  $e$.  
Then $X_e$ with the induced path metric is $\delta$-hyperbolic.
\item There exist maps ${f_e}:{X_e}{\times}[0,1]\rightarrow{X}$, such that
$f_e{|}_{{X_e}{\times}(0,1)}$ is an isometry onto the pre-image of the
interior of $e$ equipped with the path metric.
\item ${f_e}|_{{X_e}{\times}\{{0}\}}$ and 
${f_e}|_{{X_e}{\times}\{{1}\}}$ are $(K,{\epsilon})$-quasi-isometric
embeddings into $X_{v_1}$ and $X_{v_2}$ respectively for all
${v_1},{v_2}\neq{\alpha}$. When one of ${v_1},{v_2}$ is $\alpha$, this
restriction is relaxed for the corresponding inclusion of edge spaces.
${f_e}|_{{X_e}{\times}\{{0}\}}$ and 
${f_e}|_{{X_e}{\times}\{{1}\}}$ will occasionally be referred to as
$f_{v_1}$ and $f_{v_2}$ respectively.
\end{enumerate}   
\end{defn}

We shall work in the framework of pared manifolds in the sense of Thurston 
\cite{thurston-hypstr1} \cite{thurston-hypstr3}.

\begin{defn}
A {\bf pared manifold} is a pair $(M,P)$, where $P
\subset \delta M$ 
is a (possibly empty) 2-dimensional submanifold with boundary such that \\
\begin{enumerate}
\item the fundamental group of each component of $P$ injects into the
fundamental group of $M$
\item the fundamental group of each component of $P$ contains an abelian 
subgroup with finite index.
\item any cylinder $C: (S^1 \times I, \delta S^1 \times I) \rightarrow (M,P)$
with $C_\ast \colon \pi_1 (S^1 \times I )
\rightarrow \pi_1 (M)$
 injective is homotopic rel. boundary to $P$.
\item $P$ contains every component of $\delta M$ which has an abelian subgroup
of finite index.
\end{enumerate}
\end{defn}

The terminology is meant to suggest that certain parts of the skin of $M$ have been
pared off to form parabolic cusps in hyperbolic structures for
$M$. $H(M,P)$ will denote the set of hyperbolic structures on
$(M,P)$. (Note that this means that the elements of $P$ 
and the elements of $P$ alone are taken to
cusps.) 

\smallskip

{\bf Definition:} A pared manifold $(M,P)$ is said to have {\bf
  incompressible boundary} 
if each component of $\partial_0 M = \partial M \setminus P$ is
incompressible in $M$. 

Further, $(M,P)$ is said to have {\bf
  p-incompressible boundary} if \\
\begin{enumerate}
\item it has incompressible boundary \\
\item  if some curve $\sigma$ on a component of $\partial_0 M$ is freely 
homotopic in $M$ to a curve $\alpha$ on a component of $P$, then
  $\sigma$ is homotopic to $\alpha$ in $\partial M$.  \\
\end{enumerate}

\smallskip

$P_0, P_1$ will denote the components of $P$ whose fundamental groups
are virtually $Z, (Z+Z)$ respectively. The adjective `virtually' shall
  sometimes be omitted and we shall refer to the components of $P_0$
  (resp. $P_1$) as $Z$-cusps (resp. ($Z+Z$)-cusps). 

\subsection{3 Manifold as a Tree of Spaces}

The {\bf convex core} of a hyperbolic 3-manifold $N^h$ 
is the smallest convex submanifold $C({N^h}) \subset N^h$ for which inclusion 
is a homotopy equivalence.

If an $\epsilon$ neighborhood of
$C({N^h})$ has finite volume, $N^h$ is said to be
{\bf geometrically finite}. 

There exists a compact 3-dimensional submanifold
$M_{cc} \subset N^h$, the {\bf compact core} or {\bf Scott core}
\cite{scott-cc} whose inclusion is 
a homotopy equivalence. $M_{cc}$ can be thought of as $C({N^h})$ minus cusps for geometrically finite ${N^h}$.
$N^h$ minus cusps will be denoted by $N$ and $C({N^h})$ minus cusps
will be denoted by $C(N)$. The ends of $N$ are in one-to-one
correspondence with the components of 
$(N - M_{cc})$ or, equivalently, the components of
$\partial_0{M}$. 

We say that an
end of $N$ is {\bf geometrically finite} if it has a neighborhood missing
$C(N)$. 

{\bf Note:} The notion of ends here is slightly non-standard, as we do
not want to regard a cusp as an end.

An end $E$ 
of $N$ is {\bf simply degenerate} if it has a neighborhood
homeomorphic to ${S_0}{\times}{\Bbb{R}}$, where $S_0$ is the corresponding
component of $\partial_{0}{M}$, and if there is a sequence of pleated
surfaces (with cusps removed)
homotopic in this neighborhood to the inclusion of $S_0$, and exiting every
compact set. Let $S^h$ denote a hyperbolic surface of finite volume, from which $S_0$ is obtained by excising cusps.
We shall refer to $E^h = {S^h}{\times}{\Bbb{R}}$ (respecting the parametrization of $E$) as an {\it end of $N^h$}. Note
that we may think of
$E^h$ as obtained from $E$ by adjoining
"half" a $Z$-cusp.

$N$ is called {\bf geometrically tame} if all of its ends are
either geometrically finite or simply degenerate. 
Note that $N^h$ and the interior of  $N$ are homeomorphic
to the interior of $M$. For a more detailed discussion of pleated surfaces 
and geometrically tame ends, see \cite{Thurstonnotes} or
\cite{minsky-top}.

A hyperbolic structure on $(M,P)$ is a complete hyperbolic metric
on the interior of $M$ which takes precisely 
the elements of $P$  to cusps.
A manifold
$N^h$ will be said to be {\bf adapted to a pared manifold $(M,P)$} if $N^h$ corresponds
to such a hyperbolic structure on $(M,P)$.

{\bf Note:} Since the flaring ends of
$N$ contribute nothing to our discussion, we shall
(abusing notation somewhat) regard $\bf {N = C(N)}$, and refer
to $N$ as a hyperbolic manifold, though it should really
be called a convex submanifold (minus cusps) homotopy equivalent to
the big manifold (minus cusps).

A manifold $M$ has {\bf bounded geometry} if on a complement of the cusps, the
 injectivity radius
of the manifold is bounded below by some number $\epsilon $ greater
 than $0$. 
Equivalently, all closed geodesics have length greater than
 $\epsilon$.

We want to first show that the universal cover of $N^h$ minus $Z$-cusps is
quasi-isometric to a tree of hyperbolic metric spaces with possibly
exceptional vertex corresponding to the core.

Let $E^h$ be a simply degenerate end of $N^h$. Then $E^h$ is homeomorphic to
$S^h{\times}[0,{\infty})$ for some  surface $S^h$ of negative Euler
  characteristic.
Cutting off a neighborhood of the cusps of $S^h$ we get a surface with boundary
denoted as $S$. Let $E$ denote $E^h$ minus a neighborhood of the
 $Z$-cusps. We assume that each $Z$-cusp has the standard form coming from
a quotient of a horoball in ${\Bbb{H}}^3$ by $Z$. Also, we shall take
our pleated 
surfaces to be such that the pair $(S,cusps)$ is mapped to the pair
$(E,cusps)$ for each 
pleated $S^h$. We shall now show that each $\widetilde{E}$ is 
quasi-isometric to a ray of hyperbolic metric spaces satisfying the q-i
embedded condition. In \cite{mitra-trees} we had shown this for manifolds
without cusps. Each edge and vertex space will be a copy of 
$\widetilde{S}$ and the edge to vertex space inclusions shall be
quasi-isometries. Note that each $\widetilde{S}$ can be thought of as 
a copy of ${\Bbb{H}}^2$ minus an equivariant family of horodisks.

\begin{lemma}
\cite{Thurstonnotes}
There exists $D_1 > 0$ such that for all $x \in {E}$,
 there exists a pleated 
surface $g : (S^h,{\sigma}) \rightarrow E^h$ with 
$g(S){\cap}{B_{D_1}}(x) \neq \emptyset$. Also $g$ maps $(S, cusps)$
to $(E, cusps)$.
\label{closepleated}
\end{lemma}

The following Lemma
 follows easily from the fact that ${inj}_N{(x)} > \epsilon_0$:

\begin{lemma} \cite{bonahon-bouts},\cite{Thurstonnotes}
There exists $D_2 > 0$ such that if $g : (S^h,{\sigma}) \rightarrow N^h$ is a
pleated surface, then the diameter of the image of $S$ is bounded,
i.e.
$dia(g(S)) < D_2$. 
\label{diameter}
\end{lemma}

The following Lemma due to Thurston (Theorems 9.2 and 9.6.1 of
\cite{Thurstonnotes}) and Minsky \cite{minsky-top} follows from compactness
of pleated surfaces.

\begin{lemma}
\cite{minsky-top}
Fix $S^h$ and $\epsilon > 0$. Given $a > 0$ there exists $b > 0$ such that if
$g : (S^h,{\sigma})\rightarrow{E^h}$
and $h : (S^h,{\rho})\rightarrow{E^h}$ are homotopic pleated surfaces which
are isomorphisms on $\pi_1$ and $E^h$ is of bounded geometry,
then\\
\begin{center}
${d_E}(g(S),h(S)) 
\leq a \Rightarrow {d_{Teich}}({\sigma},{\rho}) \leq b$,
\end{center}
    where $d_{Teich}$ denotes Teichmuller distance.
\label{pleatedcpt}
\end{lemma}

In Lemma \ref{pleatedcpt} $d_E$ denotes the path-metric
on $E$ inherited from $N^h$. More precisely, the complete hyperbolic metric on $N^h$ gives rise to a path-metric 
on $N^h$. $E \subset N \subset N^h$ inherits a path-metric
$d_E$ given by\\
$d_E(x,y) = $ inf $\{ l(\sigma ) : \sigma $ is a path in $E$ joining $x, y \}$.

In \cite{minsky-top} a special
case of Lemma \ref{pleatedcpt} is proven for
closed surfaces. However, the main ingredient, 
a Theorem due to
Thurston is stated and proven in \cite{Thurstonnotes} (Theorems 9.2
and 9.6.1 - {\it 'algebraic limit is geometric limit'}) for finite area
surfaces. The arguments given by Minsky to prove the above Lemma from
Thurston's Theorems (Lemma 4.5, Corollary 4.6 and Lemma 4.7 of
\cite{minsky-top}) go through with very little change for surfaces
of finite area. In \cite{bowditch-ct}, Bowditch gives an alternate
approach to this using (quite general) Gromov-Hausdorff limit
arguments.

Note that in the above Lemma, pleated surfaces are not assumed to be embedded.
This is because   immersed pleated surfaces with a uniform lower
bound on  injectivity
radius are uniformly quasi-isometric to the corresponding
Riemann surfaces.

\smallskip

{\bf Construction of equispaced pleated surfaces exiting the end}\\
We next construct a sequence of {\bf equispaced} pleated surfaces
$S^h(i)\subset E^h$ exiting the end. 
 Assume that ${S^h(0)},{\cdots},{S^h(n)}$
have been constructed such that:
\begin{enumerate}
\item  $S(i),cusps$ is mapped to $E, cusps$
\item If $E(i)$ is the component of $E{\setminus}{S(i)}$ for which
$E(i)$ is non-compact, then
$S(i+1) \subset E(i)$.
\item Hausdorff distance between $S(i)$ and $S(i+1)$ is bounded above by
$3({D_1}+{D_2})$.
\item ${d_E}({S(i)},{S(i+1)}) \geq D_1 + D_2$.
\item From Lemma \ref{pleatedcpt} and condition (2)
above there exists $D_3$ depending on $D_1$, $D_2$ and $S$ such that
$d_{Teich}({S(i)},{S({i+1})}) \leq D_3$
\end{enumerate}

Next choose $x \in E(n)$, such that ${d_E}(x,{S_n}) = 2({D_1}+{D_2})$. 
Then
by Lemma \ref{closepleated}, there exists a pleated surface
$g : (S^h,{\tau}) \rightarrow E^h$ such that 
${d_E}(x,{g(S)}) \leq D_1$. Let ${S^h(n+1)} = g(S^h)$. Then by the triangle
inequality and Lemma \ref{diameter}, if $p\in{S(n)}$ and 
$q\in{S(n+1)}$,
\begin{center}
${D_1} + {D_2} \leq {d_E}(p,q) \leq 3({D_1} + {D_2})$.
\end{center}

This allows us to continue inductively. 
$S(i)$ corresponds to a point $x_i$ of $Teich(S)$. Joining the $x_i$'s in
order, one gets a Lipschitz path $\sigma$
in $Teich(S)$. 

\smallskip

\begin{defn}
 A sequence of pleated surfaces satisfying conditions
(1-5) above will be called an {\bf  equispaced sequence of pleated surfaces}.
The corresponding sequence of $S(i) \subset E$ will be called 
an {\bf equispaced sequence of  truncated pleated surfaces}.
\end{defn}

Since all
$S^h(i)$'s have bounded geometry away from cusps, they
all lie in the thick part of Teichmuller space. After
quotienting by the mapping class group, their images 
lie in a compact subset of the moduli space. Hence, by
 acting on  $S^h(i)$ by some uniformly quasi-conformal
map $\psi_i$, we may assume that $\psi_i (S^h(i))$ gives
rise to a 
fixed point $S^h$
in moduli space after
quotienting by the mapping class group. Then $\psi_{i-1}
\circ \psi_{i}^{-1} : (S^h(i)) \rightarrow (S^h(i-1))$ 
is $C$-quasiconformal for some fixed $C$. We may assume that
$\psi_i 
\circ \psi_{i-1}^{-1} : (S(i-1)) \rightarrow (S(i))$
is a bijective map by excising cusps appropriately.
Also, assume that all these excised surfaces correspond
to a fixed $S$ obtained from $S^h$ by excising cusps.

\begin{defn}
 The {\bf universal curve} over 
$X{\subset}Teich(S^h)$
is a fiber bundle over $X$ whose fiber over $x \in{X}$ is
the Riemann surface corresponding to $x$. (Topologically
this is $X{\times}S^h$.)
\end{defn}

Assuming that the $Z$-cusps are
 invariant under Teichmuller maps,
we may assume that there is an induced {\bf universal curve
of truncated hyperbolic surfaces} obtained by 
excising cusps.

Each $S(i)$ being compact (with or without boundary), $\widetilde{S(i)}$ is a hyperbolic metric space. We want to
regard the universal cover $\widetilde{E}$ of $E$ as being quasi-isometric to a  ray $T$ 
of hyperbolic metric spaces. To this end, we construct 
a quasi-isometric model of $\widetilde{E}$. Let
\begin{enumerate}
\item $T = [0,{\infty})$\\
\item vertex
set ${\mathcal{V}} = \{ n : n \in {\Bbb{N}}{\cup} \{ 0\} \}$\\
\item 
edge set 
${\mathcal{E}} = \{ [{n-1},n]: n \in {\Bbb{N}}  \}$ \\
\item 
${X_n} = \Gamma = X_{[n-1,n]}$, where $\Gamma$ is a Cayley
graph of $\pi_1(S)$ with some fixed generating set. 
\item There exist $K, \epsilon$ and a map $\eta_n$
such that $\eta_n: \Gamma \rightarrow \widetilde{S(n)}$
is a $(K, \epsilon )$ quasi-isometry for all $i$. Let 
$\eta_n^{-1}$ denote its quasi-isometric inverse.
\item The qi-embeddings
of edge sets into vertex sets are given by:\\
$\bullet$ $\phi_n : X_{[n-1,n]} \rightarrow X_n$ is the identity map on $\Gamma$
 \\
$\bullet$ $\phi_{n-1} : X_{[n-1,n]} \rightarrow X_{n-1}$
is the change of marking
 induced by sending $  \psi_{n-1}^{-1} (S)$ to
$ \psi_{n}^{-1} (S)$.\\
\end{enumerate}

By Lemma \ref{pleatedcpt} and the fact that $\sigma$
is a Lipschitz path in Teichmuller space, 
this tree of hyperbolic metric spaces satisfies the quasi-isometrically
embedded condition. In fact, we get more. Each $S(i)$ corresponds, via $S^h(i)$ to a point
$x_i$ of $Teich(S)$. Joining the $x_i$'s in
order, one gets the Lipschitz path $\sigma$ obtained above
in $Teich(S)$.  Mapping the fiber over $x_i$ to {\it an embedded incompressible
surface} lying in a (uniformly) bounded neighborhood of the corresponding pleated
surface and extending over product regions (using a metric product structure), we get a homeomorphism
between the model and $E$. Further, since the Teichmuller distance between $S^h(i)$
and $S^h(i+1)$ is uniformly bounded above, the metric product is uniformly bi-Lipschitz to
the region trapped between them in $E$. Pasting these homeomorphisms together and lifting
to the universal cover, we get

\begin{lemma}
If $E^h$ is
 a simply degenerate  end of a hyperbolic 3 manifold $N^h$
with bounded geometry, then there is a sequence of equispaced 
pleated surfaces exiting $E^h$ and hence a sequence of truncated
pleated surfaces exiting $\widetilde{E}$. Further, $\widetilde{E}$ is
quasi-isometric to a ray of hyperbolic metric spaces satisfying the 
q.i. embedded condition.
\label{equispaced}
\end{lemma}

This Lemma allows us to pass between $E$ and its quasi-isometric model,
the ray of hyperbolic metric spaces satisfying the
qi-embedded condition.

$Z$-cusps in $N$ correspond to $Z$-cusps in the boundary
 components of $\partial_0 M$. But this is not true
 for $(Z+Z)$-cusps. Recall that $P_0$ (resp. $P_1$) denotes the
 components of $P$ whose fundamental group is virtually $Z$
 (resp. $(Z+Z)$). Let 
 ${N_0} = (N \cup (Z+Z) cusps)$. 
 We shall now describe the universal cover $\widetilde{{N_0}}$ as a tree of
hyperbolic metric spaces with possibly exceptional vertex.
It is at this stage 
that we need to assume that $N^h$ is adapted to a pared manifold
$(M,P)$ with incompressible boundary. Recall
 that the incompressibility 
of the boundary $\partial_0{M}$ of a pared manifold does not require
that $\partial{M}$ be incompressible, but only that
the components of $\partial_0{M} = \partial{M} \setminus P$ be
 incompressible.

We give $(M,P_{1})$ a geometrically finite structure (with no
extra parabolics as per definition of a hyperbolic structure
adapted to $(M,P_1)$) and
denote the {\bf convex core} of this geometrically finite
manifold as $M_0$. Note that $M_0$ has no $Z$-cusps corresponding to
$P_0$ but continues to have $(Z+Z)$ cusps corresponding to
$P_1$. Further, since all $Z$-cusps have been excised in $N_0$, all
the $(Z+Z)$-cusps in $N_0$ are retained in $M_0$. 

Let ${E(1)}, {E(2)},\cdots,{E(k)}$ denote the simply 
degenerate
ends of $(N_0 - cusps)$. $M_0$ is homeomorphic by a homeomorphism that
is a quasi-isometry to (the closure of) $N_0\setminus\bigcup_i{E(i)}$.
We identify $M_0$ with its image under this map 
and denote $M_0 \cap {E(i)} = {F(i)}$, where $F(i)$ is both  
{\bf an embedded surface} in
 $E(i)$ cutting off the end, and an incompressible
 boundary component of the pared manifold
$(M_0,P_1)$. Note that $F(i)$ need not be a truncated pleated surface; but we have constructed $\widetilde{E}$
as a ray of hyperbolic metric spaces, and hence we only need
$\widetilde{F(i)}$ to be quasi-isometric to $\Gamma$, 
a Cayley graph of the 
fundamental group $\pi_1(S)$.

{\bf Remark:} That there exists such a geometrically finite hyperbolic manifold
homeomorphic to $N$ is part of Thurston's monster theorem. See 
\cite{ctm-amen} \cite{ctm-itn}
 for a different proof of the fact. Also, the limit set of
a geometrically finite manifold is locally connected \cite{and-mask}.
This shall be of use later.

{\bf Summary of Notation:} We summarize here the notation introduced
so far:\\
$\bullet$ $(M,P)$ - pared manifold \\
$\bullet$ ${P_0} \subset P$ - components of $P$ whose fundamental
group is virtually $Z$.\\
$\bullet$ ${P_1} \subset P$ - components of $P$ whose fundamental
group is virtually $(Z+Z)$.\\
$\bullet$ $N^h \in H(M,P)$. Since flaring ends of $N^h$ do not
contribute anything to the discussion, we identify $N^h$ with its
convex core $C({N^h})$.\\
$\bullet$ $N = (N^h - cusps)$\\
$\bullet$ $N_0 = N$ with $(Z+Z)$-cusps (corresponding to $P_1$)
adjoined.\\ 
$\bullet$ $M_0$ - geometrically finite structure on $(M,P_1)$.\\
$\bullet$ $E(i)$ for $i = 1 \cdots k$ denote the ends of
$N$. (Alternately, $E(i) = (E^h(i) - cusps)$). \\
$\bullet$ $M_0$ is identified with its quasi-isometric image and so is
thought of as a subset of $N_0$ \\
$\bullet$ $N_0 = M_0 \cup \bigcup_i E(i)$\\
$\bullet$ $F(i) = M_0 \cap E(i)$ \\

\begin{lemma}
 $\widetilde{N_0}$ is quasi-isometric to a tree
(T) of hyperbolic metric spaces satisfying the qi-embedded condition
with possibly exceptional vertex $\alpha$ corresponding to
$\widetilde{M_0} \subset \widetilde{N_0}$.
\label{hyptree}
\end{lemma}

{\bf Proof:}
Note that $\widetilde{M_0} \subset \widetilde{N_0}$ is the universal
cover of the convex core of a geometrically finite manifold and hence
 is a hyperbolic
metric space. Let $\widetilde{F(i)} \subset \widetilde{N_0}$ represent a lift
of $F(i)$ 
to $\widetilde{N_0}$. Then, $\widetilde{F(i)}$, being quasi-isometric 
to the fundamental group of a compact surface
(with or without boundary)
 is a word-hyperbolic metric space. If $\widetilde{E(i)}$ is
a lift of $E(i)$ containing $\widetilde{F(i)}$ then from Lemma
\ref{equispaced},
$\widetilde{E(i)}$ is a ray of hyperbolic metric spaces satisfying
the q.i. embedded condition. Since there are only
finitely many ends $E_i$, we can thus regard $X = \widetilde{N_0}$ as
a tree ($T$) of hyperbolic metric spaces such that \\
$\bullet$ $\alpha$ is the root of $T$. $X_\alpha = \widetilde{M_0}$ is a
  hyperbolic metric space.\\
$\bullet$ $T$ consists of a (finite or infinite) collection of rays
  emanating from $\alpha$.\\
$\bullet$ Each copy of $\widetilde{E(i)}$ in $\widetilde{N_0}$ is
  quasi-isometric to a ray ($T_i$) of hyperbolic metric spaces
  satisfying the q.i. embedded condition (Lemma \ref{equispaced}).\\
$\bullet$ $N_0 = M_0 \cup \bigcup_i E(i)$ \\
$\bullet$ $F(i) = M_0 \cap E(i)$ \\
$\bullet$ As of now no restrictions are imposed on the inclusion of each
  $\widetilde{F_i}$ into $\widetilde{M_0}$ \\

These are precisely the defining conditions of a tree of hyperbolic
metric spaces with possibly exceptional vertex satisfying the
q.i. embedded condition. $\Box$

\smallskip

{\bf Note:} $T$ has a root vertex $\alpha$ which is possibly
exceptional. The rest of $T$ consists of
 a number of rays emanating from $\alpha$.

\subsection{A Topological Property of Pared Manifolds}

Before we enter into the construction of quasiconvex sets, we shall
describe a basic topological property of pared manifolds.

\begin{lemma}Let $(M,P)$ be a pared manifold with incompressible
boundary. Then \\
$\bullet$ No annulus component of $P$ is freely homotopic to a 
 curve on a torus component or on another annulus component of $P$ \\
$\bullet$ If two curves (which are not non-trivial powers of any other curves)
on $\partial_0 M = \partial M \setminus P$
are freely homotopic to curves on the same torus component 
of $P$, then they are in fact freely homotopic to the same curve on a
torus component of $P$ and hence to each other. \\
\label{pared}
\end{lemma}

{\bf Proof:} {\it Statement 1:}  Let $\Delta$ be a boundary
torus.  If possible, 
let $A$ be an annulus in the boundary of $M$ such that its
core curve is freely homotopic to a curve $\sigma$ on $\Delta$. The complement
of a small neighborhood of $\sigma$ in
$\Delta$ is an annulus $A_1$. Connecting the boundary curves  of $A_1$
to those of $A$ by the free homotopy  we get an immersed annulus with
boundary on the boundary of $M$ , but not homotopic to $\partial{M}$. This
proves that the core curve of $A$ cannot be homotopic to the core curve
of an annular component of $P$ as this would contradict 
condition 3 of the definition of a pared manifold. 
If two core curves of annuli components of $P$ are homotopic, we get
a new annulus interpolating between these curves, again contradicting
Condition 3 of the definition of a pared manifold. This proves the
first part of the lemma.

\smallskip

\noindent {\it Statement 2:} Let $M_\Delta$ denote
the cover of $M$
corresponding to $\pi_1{({\Delta})}$. 
Let
$\sigma_1$ and $\sigma_2$ be the curves
 homotopic to curves $\alpha_1$
and $\alpha_2$ on $\Delta$. Let $B$ denote the annulus between
$\sigma_1$ and $\alpha_1$, then the intersection number between $\alpha_1$
and $\alpha_2$ in $\Delta$ is the same as the intersection number  between
 $\alpha_2$ and $B$ in $M_\Delta$ (see for instance \cite{Thurstonnotes} or \cite{bonahon-bouts}). 
Since $\alpha_2$, $\sigma_2$ are homotopic in 
$M$, this homotopy lifts to $M_\Delta$.
Hence, the intersection number  between
 $\alpha_2$ and $B$ in $M_\Delta$ is the same as the
intersection number between $\sigma_2$ and $B$
in $M_\Delta$ (by homotopy invariance of
intersection numbers). 
If $\sigma_1$ and $\sigma_2$
are on different components then this latter number is zero. 

Else, they
belong to the same component $K$ of $\partial_0{(M)}$. $\pi_1 (K)$ is
either free non-abelian or else $K$ is a closed surface of genus
greater than $1$. Since the fundamental
group of $K$ injects into that of $M$, and since elements represented
by $\sigma_1$, $\sigma_2$ commute, therefore $\sigma_1$ and $\sigma_2$
must be powers of the same curve. But the hypothesis says that the curves
$\sigma_1$ and $\sigma_2$ are not powers of any other curves. Hence the two
denote the same curve. 

In either case, $\sigma_1$ and $\sigma_2$ are freely homotopic to the
same curve on $\Delta$ and hence to each other.$\Box$

\section{Construction of q.i. embedded sets }

From Lemma \ref{hyptree} we have that  $\widetilde{N_0}$ is quasi-isometric
 to a tree
(T) of hyperbolic metric spaces satisfying the qi-embedded condition
with possibly exceptional vertex $\alpha$ corresponding to
$\widetilde{M_0} \subset \widetilde{N_0}$. In fact
$T$ has a root vertex $\alpha$ and consists of 
 a number of distinct rays emanating from $\alpha$.
Let $X$ denote the tree of spaces, $X_v$ denote vertex space
 corresponding to vertex $v$, $X_e$ denote edge space corresponding to
 edge $e$.

For convenience of exposition, we shall sometimes need to 
modify $X$, ${X_v}$, ${X_e}$ by
 quasi-isometric perturbations and regard them as graphs. Given a
 geodesically complete metric 
space $(Z,d)$ of bounded geometry, choose a maximal disjoint collection  
$\{{B_1}({z^{\prime}})\}$ of disjoint 1-balls. Then by maximality,
for all $z \in Z$ 
there exist $z^\prime$ in the collection such that $d(z,{z^\prime}) < 2$.
Construct a graph $Z_1$ with vertex set $\{{z'}\}$ and edge set consisting
of distinct vertices ${z_{a}}$, ${z_{b}}$ 
such that $d({z_{a}},{z_{b}}) \leq 4 $. Then $Z_1$ equipped with
the path-metric is quasi-isometric to $(Z,d)$
(see for instance \cite{gromov-greenbk}). Henceforth we shall
move back and forth between descriptions of spaces as Riemannian
manifolds and as graphs, assuming that there is a quasi-isometry
between them. The quasi-isometry will usually be suppressed.

Let $v$ be a vertex of $T$. 
Let $v_{-} \neq\alpha $ be the penultimate vertex on the geodesic edge path
from $\alpha$ to $v$. Let $e_-$ denote the directed edge from ${v_{-}}$
to $v$.  Define a quasi-isometry
${\phi_v} : {f_{e_{-}}}({X_{e_{-}}}{\times}\{{0}\}) \rightarrow
{f_{e_{-}}}({X_{e_{-}}}{\times}\{{1}\})$ as follows: \\
If $p{\in}{f_{e_{-}}}({X_{e_-}}{\times}\{0\}){\subset}{X_{v_{-}}}$, choose
$x\in{X_{e_{-}}}$ such that $p={f_{e_{-}}}(x{\times}\{0\})$ and define
\begin{center}
${\phi_v}(p) =
{f_{e_{-}}}({x}{\times}\{{1}\})$.
\end{center}

Note that in the above definition, $x$ is arbitrarily chosen from a
set of bounded 
diameter.

Let $\mu$  be a geodesic in 
$X_{v_{-}}$, joining $a, b
 \in {f_{e_{-}}}({X_{e_{-}}}{\times}\{{0}\}) $. ${\Phi}_{v}({\mu})$
will denote a geodesic in $X_v$ joining $\phi_v{(a)}$ and $\phi_v{(b)}$.
For our purposes since all the edge and vertex spaces on a ray, (apart
from $X_\alpha$) are identical, we might as well identify 
 ${f_{e_{-}}}({X_{e_-}}{\times}\{0\}){\subset}{X_{v_{-}}}$ with 
$X_{v_-}$ itself. Similarly we may identify
 ${f_{e_{-}}}({X_{e_-}}{\times}\{1\}){\subset}{X_{v}}$ with $X_v$
itself. Thus, $\phi_v$ 
is regarded as a quasi-isometry from $X_{v_{-}}$ to $X_v$ for $v, v_-
\neq \alpha$. We define a corresponding map $\Phi_v$ from geodesics in
$X_{v_-}$ to geodesics in $X_v$ by taking a geodesic joining $a, b \in
X_{v_-}$ to one joining $\phi_v{(a)},\phi_v{(b)} \in X_v$. 
From Lemma \ref{equispaced}, there exist $k, \epsilon > 0$ such that
for all $v$, $\phi_v$ is a $(k, \epsilon)$ - quasi-isometry.

\smallskip

In Section 4.1, we shall construct a {\bf hyperbolic ladder-like} set
$B_\lambda$ containing $\lambda$. In Section 4.2, we shall construct a
retract $\Pi_\lambda : \widetilde{N_0} \rightarrow B_\lambda$. In
Section 4.3, we shall prove that $\Pi_\lambda$ {\it is} a retract,
i.e. it fixes $B_\lambda$ and does not stretch distances much. This
will show that $B_\lambda$ is quasi-isometrically embedded in
$\widetilde{N_0}$.

\subsection{Construction of the hyperbolic ladder-like set $B_\lambda$}

\medskip

The quasi-isometrically embedded set $B_\lambda$ that we intend to
construct
will contain the images of a geodesic under such quasi-isometries.
 Suppose $F(i)$ cuts off the end $E(i)$ and
$\mu\subset\widetilde{F(i)}\subset \widetilde{E(i)}$. Then denote the union of
the images of $\mu$ under the quasi-isometries taking it to the
different vertex spaces as $B({\mu})$. That is to say, if $v_0,
v_1,\cdots $ denote the vertices of the ray exiting the end, then
the union of $\mu$, $\Phi_{v_1}({\mu})$,
$\Phi_{v_2}\circ\Phi_{v_1}({\mu})$,
$\Phi_{v_3}\circ\Phi_{v_2}\circ\Phi_{v_1}({\mu})$, etc. is denoted
by $B({\mu})$.

\medskip

Recall that $M_{0}$ denotes a convex geometrically finite hyperbolic manifold
with hyperbolic structure adapted to the pared manifold $(M,P_1)$,
where
$P_1 \subset P$ denotes the set of $(Z+Z)$-cusps. Also
$\widetilde{M_{0}}$
the universal cover of $M_{0}$ is identified with $\widetilde{N_0}
- \bigcup_i E(i)$. $M_{cc}$ will denote $M_{0}$ minus a
neighborhood of the $(Z+Z)$-cusps. Note that $M_{cc}$ can be thought
of as the Scott core of the manifold. $M_{gf}$ will
denote a geometrically finite hyperbolic structure 
 adapted to the pared manifold $(M,P)$. The difference between
 $M_{gf}$ and $M_0$ is that while $M_{gf}$ is adapted to the pared
 manifold $(M,P)$, $M_0$ is adapted to the pair $(M,P_1)$. Thus
 $M_{gf}$ has incompressible boundary as a pared manifold, but the
 same may not be true of $M_0$. One may also think of $M_{gf}$ as
 $M_0$ with $Z$-cusps (corresponding to $P_0$) adjoined.

Fix a set of neighborhoods of the cusps of $M_0$, which are
sufficiently
separated from each other. Recall that
 $M_{cc}$ denotes $M_0$ minus these cusps. Then
$\widetilde{M_{cc}} = \widetilde{M_0} \setminus \bigcup \mathcal{H}$
where $\mathcal{H}$ denotes an equivariant system of horoballs
corresponding to lifts of $(Z+Z)$-cusps. Note that $\widetilde{M_{cc}}$
is quasi-isometric to the Cayley graph of ${\pi_1}(M)$ as the quotient is compact \cite{GhH} \cite{gromov-hypgps}. 

Let $\lambda^h = [a,b]$ be a hyperbolic geodesic in $\widetilde{M_{gf}}$.
Let ${\beta}^h$ be a hyperbolic geodesic in $\widetilde{N^h}$ joining $a,
b$. Here
$\widetilde{M_{gf}}$ is identified with its image in
$\widetilde{N^h}$. We shall show that if $\lambda^h$ lies outside  a large
ball around a fixed reference point $p\in \widetilde{M_{gf}}$ then so
does ${\beta}^h \in \widetilde{N^h}$. 
Recall that $F(i)$ for $i = 1 \cdots k$ are components of $\partial_0 M
= \partial{M} - P$. We can identify each $F(i)$ with a boundary
component of $M_{0}$ so that the inclusion of $F(i)$ into $M_{cc}$
induces an injection at the level of the fundamental group. We
identify each $F(i)$ with the first {\bf truncated pleated surface}
exiting the end $E(i)$.

\medskip

{\bf Summary of Notation:} \\
$\bullet$ $M_0$ - convex geometrically finite hyperbolic structure
on $(M,P_1)$\\
$\bullet$ $M_{gf}$ -  convex geometrically finite hyperbolic
structure adapted to the pared manifold $(M,P)$ \\
$\bullet$ $M_{cc}$ - $M_0$ minus $(Z+Z)$-cusps \\
$\bullet$ $N^h \in H(M,P)$ has bounded geometry and is geometrically tame.\\
$\bullet$ $N$ - $(N^h - cusps)$\\
$\bullet$ $N_0$ - $N$ with $(Z+Z)$ cusps adjoined\\
$\bullet$ $M_{gf}$ is identified with its homeomorphic image in $ N^h $
taking cusps to cusps\\
$\bullet$ $M_{0}$ is identified with its homeomorphic image in $N^0$\\
$\bullet$ $M_{cc}$ is identified with its homeomorphic image in $N$\\

\medskip

{\underline{\bf Construction of $B_0{( \lambda )}$}}\\
Given $\lambda^h$ we shall first construct an ambient quasigeodesic
$\lambda$ in $\widetilde{M_0}$.  Since $\lambda^h$ is a hyperbolic
geodesic in  $\widetilde{M_{gf}}$ there are unique entry and exit
points for each horoball that $\lambda^h$ meets and hence unique
Euclidean geodesics joining them on the corresponding
horosphere. Replacing the segments of $\lambda^h$ lying inside
$Z$-horoballs by the corresponding Euclidean geodesics, we obtain an
ambient quasigeodesic $\lambda$ in $\widetilde{M_0}$ by Theorem
\ref{ctm}. See Figure 1 below:

\medskip

\begin{center}
\includegraphics{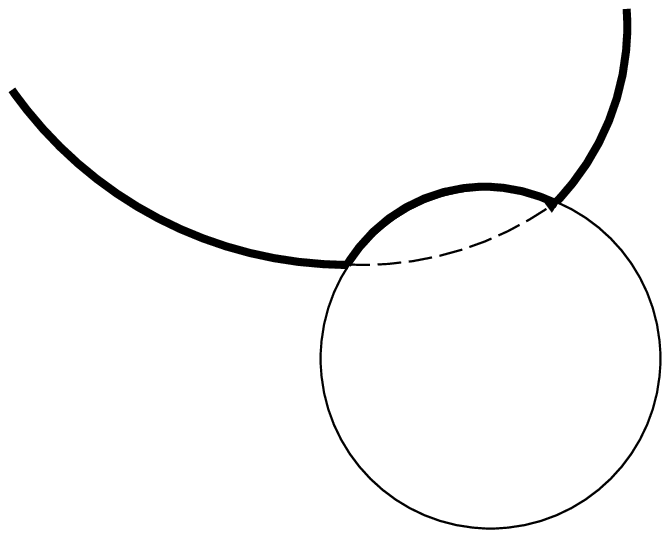}

\smallskip

\underline{Figure 1}

\end{center}

\smallskip

Again, since $\lambda$ coincides with $\lambda^h$ outside horoballs
corresponding to $Z$-cusps, there exist unique entry and exit points
of $\lambda$ into horoballs corresponding to $(Z+Z)$-cusps, and hence
again, unique Euclidean geodesics joining them on the corresponding
horosphere.  Replacing the segments of $\lambda$ lying inside
$(Z+Z)$-horoballs by the corresponding Euclidean geodesics, we obtain an
ambient quasigeodesic $\lambda_{cc}$ in $\widetilde{M_{cc}}$ by Theorem
\ref{ctm}. Each
of the Euclidean geodesics mentioned above
along with the hyperbolic geodesic
joining its end-points (and lying entirely within the horoball) bounds
a 
 totally geodesic disk. (See Figure 2 below)

\medskip

\begin{center}
\includegraphics{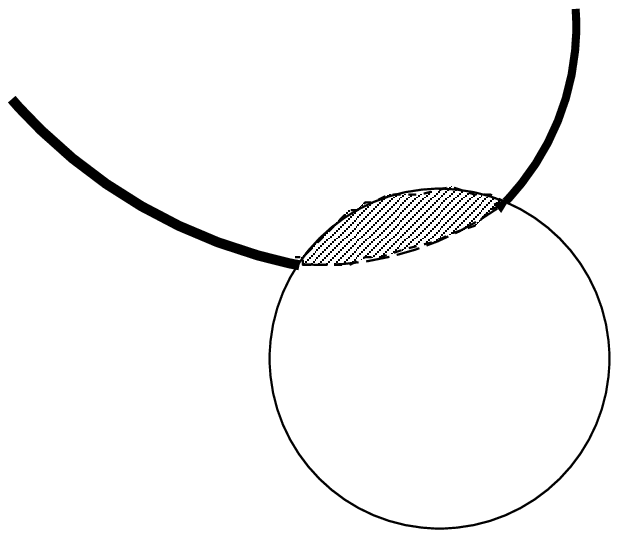}

\smallskip

\underline{Figure 2}

\end{center}

\smallskip

The union of $\lambda$ and all these totally
geodesic disks is denoted by $B_0{( \lambda )}$. There exists $C_1 >
0$ such that each $B_0{({\lambda})}$ is $C_1$-quasiconvex in
$\widetilde{M_0}$. (See for instance McMullen \cite{ctm-locconn}
Section 8.) 

$B_0^{aug} ( \lambda )$ will denote $\lambda \cup \mathcal{H} (
\lambda )$, where $\mathcal{H} ( \lambda ) \subset \mathcal{H}$
denotes the collection of horoballs in $\mathcal{H}$ that $\lambda$
meets. Again, from Theorem \ref{ctm}, $B_0^{aug} ( \lambda )$ can be
assumed to be $C_1$-quasiconvex in $\widetilde{M_0}$.

{\bf Note:} We will be using $\widetilde{N_0}$ rather than
$\widetilde{N^h}$ for the construction of $B_{\lambda}$. Hence
$\widetilde{M_0}$ 
rather than $\widetilde{M_{gf}}$  is the relevant space. Recall (Lemma
\ref{hyptree} ) that $\widetilde{N_0}$ is a tree of hyperbolic metric
spaces with possibly exceptional vertex corresponding to
$\widetilde{M_0}$ satisfying the q.i. embedded property. 

\medskip

{\underline{\bf Construction of $B_1{( \lambda )}$}}\\
Technically, (though not conceptually) 
this step is the most intricate one, and is the one new construction
that is required to handle parabolics. The values of the constants
chosen here become clear only through hindsight. The reason behind the
choices made here will become clear only while constructing the
projection in the next subsection.
As mentioned in Section 1.1, we
have to be quite careful while handling the different kinds of cusps
that arise:

\begin{enumerate}
\item {\bf $Z$-cusps} \\
\item {\bf $(Z + Z)$ cusps}, where no curve on any component of
$\partial_0 M$ is homotopic to a curve on the boundary torus
corresponding to the  $(Z + Z)$ cusp.\\
\item {\bf $Z$ and  $(Z + Z)$ cusps}, where some curve(s) on component(s) of 
$\partial_0 M$ is(are)  homotopic to a particular  curve on the boundary torus
corresponding to the  $(Z + Z)$ cusp or to some multiple of the core
curve of a $Z$-cusp. 
\end{enumerate}

Next we need to consider the parts of $\lambda$ that can have
substantial overlap with the ends $\widetilde{E}$, that is to
say those pieces of $\lambda$ that follow some $\widetilde{F(i)}$ 
 for a considerable length. There are geodesic segments on each 
 $\widetilde{F(i)}$ parallel to these pieces. The union of all these
parallel geodesic segments along with $B_0{({\lambda})}$ will be
denoted by $B_1{({\lambda})}$. Details follow.

\begin{lemma}
Let $X$ be a convex subset of ${\Bbb{H}}^n$. Let $Y, Z$ be closed
convex subsets of $X$. Let $Y_\delta = {N_\delta}(Y)$, $Z_\delta =
{N_\delta}(Z)$ be closed $\delta$-neighborhoods of $Y, Z$
respectively where $\delta$ denotes the hyperbolicity constant of the
ambient space (in this case a convex subset of ${\mathcal{H}}^n$). 
Let $\Pi$ denote nearest point projection of $X$ onto
$Y$. Then $d_X(\Pi{(Z)},{Y_\delta}{\cap}{Z_\delta})$ is
bounded in terms of $\delta$.
\label{projn-bdd}
\end{lemma}

\noindent {\bf Proof:} {\bf Case (a):}  $({Y_\delta}){\cap}({Z_\delta}) \neq
\emptyset$\\
Let $z \in Z$. We want to show that there exists $D = D( \delta
) > 0$ 
independent of $z$ such
that  $d_X(\Pi{(z)},({Y_\delta}){\cap}({Z_\delta})) < D$. Let $\Pi_0$
denote nearest point retraction of $X$ onto the convex set 
 $({Y_\delta}){\cap}({Z_\delta})$ (which, being an intersection of
convex sets, is convex.)  

As usual, $[a,b]$ denotes a geodesic joining $a, b$.

\begin{eqnarray*}
[z,{\Pi}(z)] & \subset &
N_\delta{([z,{\Pi_0}(z)]\cup{[{\Pi_0}(z),{\Pi}(z)]})}  
\end{eqnarray*}

But $[z,{\Pi_0}(z)] \subset Z_\delta$ and 
${[{\Pi_0}(z),{\Pi}(z)]} \subset Y_\delta$. Hence, 

\begin{eqnarray*}
[z,{\Pi}(z)]& \subset & {N_\delta}(Y{\cup}Z)
\end{eqnarray*}

\noindent
$\bullet$ $z \in Z \subset Z_\delta$\\
$\bullet$ $\Pi_0{(z)} \in Y_\delta{\cap}Z_\delta$\\
$\bullet$ $\Pi{(z)} \in Y_\delta$ \\

The geodesic $[z,{\Pi}(z)]$ has to cross into $Y_\delta$ at some point
$p \in {Y_\delta}\cap{Z_\delta}$. (Note that here we are using the fact
that ${Y_\delta}\cap{Z_\delta}$ is closed.) But then $d(\Pi (z),p) \leq
\delta$ (else we would be able to find another point $q$ at a distance
of less than $\delta$ from $p$ in $Y$ contradicting the definition of
${\Pi}(z)$.) This proves Case (a).

\smallskip

{\bf Case (b):} ${Y_\delta}\cap{Z_\delta} = \emptyset$ \\
Let $z_1, z_2 \in Z$ and let $y_i = \Pi{(z_i)}$ for $i = 1, 2.$ Since
local quasi-geodesics are global quasi-geodesics \cite{gromov-hypgps}
\cite{GhH}
in hyperbolic metric
spaces, there exists $D > 0$ such that if $d(y_1,y_2) \geq D$, then
$[z_1,y_1]{\cup}[y_1,y_2]{\cup}[y_2,z_2]$ is a
$2 \delta$-quasigeodesic. (To see this,
first note that  Lemma 3.1 of \cite{mitra-trees} 
ensures that $[z_1,y_1]{\cup}[y_1,y_2]{\cup}[y_2,z_2]$ is a
$C = C(2 \delta )$-quasigeodesic.That $C = 2 \delta$ comes from applying $\delta$-thinness to
triangles $(z_1, y_1, y_2)$ and $( y_1, y_2, z_2)$ 
in succession.)\\
 In this case,
$([z_1,y_1]{\cup}[y_1,y_2]{\cup}[y_2,z_2])\subset{N_{2 \delta}}[z_1,z_2]$.
But $[z_1 ,z_2] \subset Z$ as $Z$ is convex. In particular $y_1, y_2
\in N_{2\delta}(Z)$. This contradicts the assumption that 
 ${Y_\delta}\cap{Z_\delta} = \emptyset$ \\
Hence, $d(y_1,y_2) \leq D$, proving the Lemma. $\Box$

\smallskip

The above Lemma can be `quasi-fied' as follows. (The scheme of proof
being a `quasification' of the above, we omit it.)

\begin{lemma}
Given $ \delta , C > 0$, there exists $k > 0$ such that the following
holds: \\
Let $X$ be a $\delta$-hyperbolic metric space and $Y, Z$ be
$C$-quasiconvex subsets of $X$. Let $Y_k, Z_k$ denote
$k$-neighborhoods of $Y, Z$ respectively. Let $\Pi$ denote a nearest
point projection of $X$ onto $Y$. Then $d_X({\Pi}(Z), (Y_k\cap{Z_k}))$
is uniformly bounded. Hence there exists $k > 0$ such that
 $({\Pi}(Z) \subset (Y_k \cap {Z_k}))$ if the latter set is non-empty
and $\Pi (Z)$ has diameter bounded by $k > 0$ if $Y_k \cap Z_k =
\emptyset$.
\label{projn-bdd2}
\end{lemma}

Let $F$ be one of the $F(i)$'s. Either
$\widetilde{F}$ is quasiconvex in $\widetilde{{M_0}}$ (when there are
no accidental parabolics) or there exist
disjoint curves $\sigma_i$ for $i = 1 \cdots l$ (by Lemma \ref{pared}) 
homotopic to curves on  components of $P$. Therefore $\sigma_i$
correspond to parabolics. Each $\sigma_i$ lies at a bounded distance
from some Euclidean geodesic $\eta_i$ on a  torus or annular component
of $P$ in
$M_{cc}$. Each $\eta_i$ forms the boundary of
 a totally geodesic $Z$ cusp 
$\kappa_i$ (possibly 
a totally geodesic subset of a $Z+Z$-cusp). Further, the two curves $\eta_i$ and $\sigma_i$ bound
between themselves an annulus $A_i$. Let $G = F \cup \bigcup_i
({A_i}{\cup}{\kappa_i})$. Then $\widetilde{G}$ is quasiconvex in
$\widetilde{M_{0}}$. Choose a  constant $C_2 > 0$ such that each such
$\widetilde{G}$ is $C_2$-quasiconvex in $\widetilde{M_{0}}$. Recall
 that $B_0{({\lambda})}$ and  $B_0^{aug}{({\lambda})}$ are
 $C_1$-quasiconvex.

We are now in a position to start constructing  $B_1{({\lambda})}$ from
 $B_0{({\lambda})}$. 

Choosing $ Y =  B_0{({\lambda})}$ and $Z = \widetilde{G (i)}$ (one
lift of $G(i)$ is considered in turn, and choosing $k$ to be the
maximum of the corresponding finite collection of $k$'s), we obtain a
$k > 0$ from Lemma \ref{projn-bdd2} above,  so that 
 $({\Pi}(Z) \subset (Y_k\cap{Z_k}))$. Let $\widetilde{F(i)}$, $i = 1
\cdots s$ denote the different lifts of $F_i$'s that intersect a $(3k
+ 4\delta )$ neighborhood of $B_0 ( \lambda )$. Since $\lambda$ is
a finite segment, the number $s$ is finite. Choose $p_i, q_i \in
\widetilde{F(i)} \cap N_{3k + 4\delta }(B_0 ( \lambda ))$ such that $d({p_i},{q_i})$ is maximal.

Fix $D > 0$. Choose the copies of $\widetilde{F(i)}$ for which 
  $d({p_i},{q_i}) \geq D$. 

\smallskip

{\bf Note:} The number $D$ picks up significance in the {\it
  p-incompressible} case. We shall come back to this in 
  Section 5.2, when we simplify our problem under this extra assumption.

Redefining $s$ if
  necessary, we let $\widetilde{F(1)}\cdots\widetilde{F(s)}$ denote this
  collection. Let $\mu_i$ denote the
  geodesic in $\widetilde{F(i)}$ joining $p_i$ and $q_i$. 
Let $\mathcal{E}$ denote the corresponding collection of
  $\widetilde{E(i)}$'s.  Let ${\mathcal{E}}^{\prime}$ denote 
the collection of
  the remaining  $\widetilde{E(i)}$'s. 
We would like to define

\smallskip

\begin{center}

 $B_1{({\lambda})}=B_0{({\lambda})}\cup \bigcup \mu_i$. 

\end{center}

\medskip

The choice of $\mathcal{E}$, ${\mathcal{E}}^{\prime}$ is dependent on
$k$. We want each $\widetilde{E}$ to be such that \\

\smallskip

$\bullet$ Either $\widetilde{F}$ is quasiconvex in $\widetilde{M_0}$
(case where no curve on $F$ is parabolic in $\widetilde{M_0}$) .

$\bullet$ or, $\widetilde{G}$ (which is always quasiconvex in
$\widetilde{M_0}$) intersects at most one 
 $\bf{H} \in \mathcal{H} ({\lambda})$

\smallskip

\begin{lemma}
There exists $k^\prime > 0$ such that if $\widetilde{G}$ intersects
more than one horoball in ${\mathcal{H}}
({\lambda})$, then $N_{k^{\prime}} (B_0 (
\lambda )) \cap \widetilde{F}$ has non-zero diameter. 
\label{oneball}
\end{lemma}

{\bf Proof:} Suppose $\widetilde{G}$ intersects more than one horoball
$\bf{H} \in \mathcal{H} ( {\lambda} )$, say ${\bf{H}}_1$,
${\bf{H}}_2$. Let $p \in  {\bf{H}}_1 \cap \widetilde{G}$ and
 $q \in  {\bf{H}}_2 \cap \widetilde{G}$. Let $[p,q]_0$ denote the
geodesic in $\widetilde{G}$ joining $p, q$, and let $[p,q]$ denote the
geodesic in $\widetilde{M_0}$ joining $p, q$. $[p,q]_0$ lies in a
$C_1$-neighborhood of $[p,q]$ as $\widetilde{G}$ is
$C_1$-quasiconvex. Since horoballs are sufficiently separated in both
$\widetilde{G}$ and $\widetilde{M_0}$, therefore some part of
$[p,q]_0$ lies on $\widetilde{F}$. Hence $N_{C_1} ([p,q]) \cap
\widetilde{F} \neq \emptyset$ and has non-zero diameter. 

Next, there is a subsegment of $\lambda$ 
starting in  ${\bf{H}}_1$ and ending in  ${\bf{H}}_2$. Let
$\lambda^{\prime}$ be the part of this subsegment lying outside the
interior of the horoballs ${\bf{H}}_1$ and  ${\bf{H}}_2$. Then, off
horoballs, $\lambda^{\prime}$ and $[p,q]$ track each other (by Theorem
\ref{farb}), i.e. there exists $C_2 > 0$ such that $[p,q] \subset
N_{C_2} ( \lambda \cup \mathcal{H} ( \lambda ))$ and the subsegment of
$[p,q]$ lying outside  the horoballs ${\bf{H}}_1$ and  ${\bf{H}}_2$
 is contained in a $C_2$-neighborhood of $B_0 ( \lambda )$. Choosing
 $k^{\prime} = C_1 + C_2$, we are through. $\Box$

\smallskip

The next corollary follows:

\begin{cor}
Given $k > 0$, there exists $k^\prime > 0$ such that if the
$2k$-neighborhood of $\widetilde{G}$ intersects
more than one horoball in ${\mathcal{H}}
({\lambda})$, then the diameter of
 $N_{k^{\prime}} (B_0 (
\lambda )) \cap \widetilde{F}$ is non-zero. 
\label{oneball2}
\end{cor}

Here the term $2k$ occurs so as to recall the $k$ of Lemma
\ref{projn-bdd2}. 
Choose $K = 3k + 4\delta + k^\prime$, where $k$ is as in Lemma
\ref{projn-bdd2} and $k^\prime$ is as in Corollary
\ref{oneball2} above. 

We now return to our construction of $B_1( \lambda )$, or more
specifically the $\mu_i$'s that we had mentioned after Lemma
\ref{projn-bdd2}. 

 Let $\widetilde{F(i)}$, $i = 1
\cdots s$ (redefining $s$ if necessary) 
 denote the different lifts of $F(i)$'s that intersect a $K$
 neighborhood of $B_0 ( \lambda )$. Since $\lambda$ is 
a finite segment, the number $s$ is finite. Choose $p_i, q_i \in
\widetilde{F(i)}\cap N_{3k + 4\delta }(B_0 ( \lambda ))$ such that $d({p_i},{q_i})$ is maximal.

Recall that
$p_i, q_i \in
\widetilde{F(i)} \cap N_{3k + 4\delta }(B_0 ( \lambda ))$ for some copies of $\widetilde{F(i)}$.
 Choose the copies of  $\widetilde{F(i)}$ for which 
  $d({p_i},{q_i}) > 0$. Redefining $s$ once more, if
  necessary, we let $\widetilde{F(1)}\cdots\widetilde{F(s)}$ denote this
  collection. Let $\mu_i$ denote the
  geodesic in $\widetilde{F(i)}$ joining $p_i$ and $q_i$. 
Let $\mathcal{E}$ denote the corresponding collection of
  $\widetilde{E(i)}$'s.  Let ${\mathcal{E}}^{\prime}$ denote 
the collection of
  the remaining  $\widetilde{E(i)}$'s. 
We are finally in a position to  define

\smallskip

\begin{center}

 $B_1{({\lambda})}=B_0{({\lambda})}\cup \bigcup \mu_i$. 

\end{center}
See Figure 3 below. $\lambda$ lies on $\widetilde{M_0}$. `Parallel'
segments $\mu_1, \cdots \mu_k$ are constructed lying on
$\widetilde{F(1)}, \cdots \widetilde{F(k)}$.

\medskip

\begin{center}
\includegraphics[height=3.5cm]{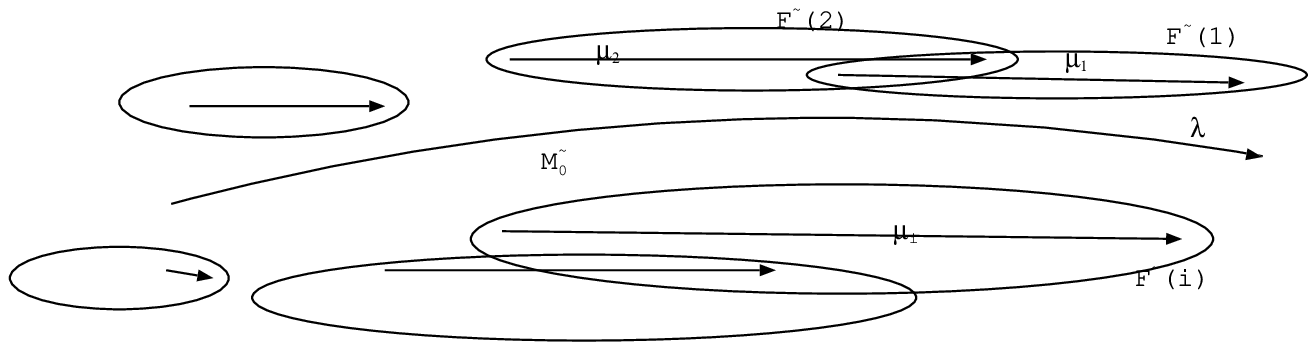}

\smallskip

\underline{Figure 3}
\end{center}

\medskip

{\underline{\bf Construction of $B{( \mu )}$}}\\
Recall from Lemma \ref{equispaced} that each $\widetilde{E}$ is
quasi-isometric to a ray $[0,{\infty})$ of hyperbolic metric spaces
  with integer points corresponding to vertices and $[n-1,n]$ (for
  $n{\in}\Bbb{N}$ ) corresponding to edges. Fix a particular $\widetilde{F}$
  cutting off the end $\widetilde{E}$ and a geodesic segment
  $\mu\subset\widetilde{F}$. Let $X_i$ denote the vertex spaces for $i
  = 0, 1, \cdots$. Let $\phi_i$ denote the quasi-isometry from $X_{i-1}$
  to $X_i$ and $\Phi_i$ denote the corresponding map from geodesic
  segments in  $X_{i-1}$
  to those in $X_i$. Define 

\begin{center}
$\mu^i =
  \Phi_i\circ\cdots\circ\Phi_1{({\mu})}$\\
 and $B{({\mu})} = \bigcup\mu^i$
\end{center}
See Figure 4 below.

\medskip

\begin{center}
\includegraphics{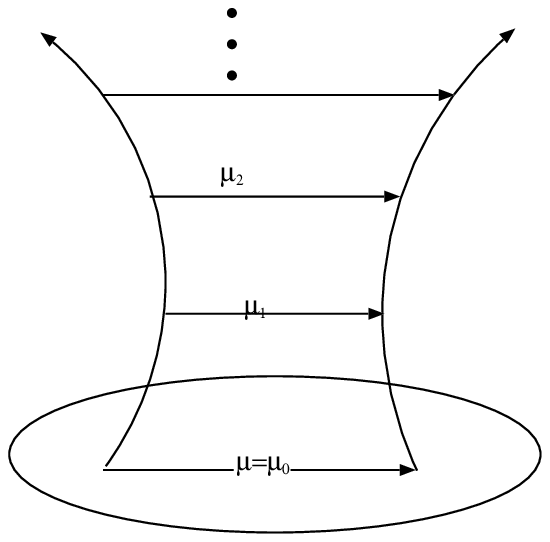}

\smallskip

\underline{Figure 4: {\it Hyperbolic ladder-like Set} }
\end{center}

\smallskip

{\underline{\bf Definition of $B_\lambda$}}

\noindent Finally define

\begin{center}
$B_\lambda = B_1{({\lambda})}\cup\bigcup_{i=1\cdots{s}}{B{({\mu_i})}}$
\end{center}

\subsection{Definition  of $\Pi_\lambda$}

\medskip

Recall ${\mathcal{E}} = \{ \widetilde{E(1)} \cdots \widetilde{E(s)}\}$ and
 ${\mathcal{E}}' = \{ {\widetilde{E(s+1)}}, {\widetilde{E(s+2)}},
 \cdots \}$. 
We next want to show that $B_\lambda$ is quasi-isometrically
embedded. To prove this, we shall construct a retraction
$\Pi_\lambda$. Let $X_{ij}$ denote the vertex spaces in 
$\widetilde{E(i)}\in\mathcal{E}$. Below,
 we shall have need
to replace $\widetilde{E(i)}$ by the union
of vertex spaces $\bigcup_jX_{ij}$ (fixing $i$
and letting $j$ vary from $0$ to $\infty$)
with the understanding that the latter has
the  metric induced from
 $\widetilde{E(i)}$. The difference between $\widetilde{E(i)}$ and
$\bigcup_jX_{ij}$ is just that the edge spaces
of $\widetilde{E(i)}$
are not explicitly present in the latter.

Let $\phi_{ij}$ denote the
 quasi-isometry from $X_{i,j-1}$ to $X_{ij}$ and let $\Phi_{ij}$
 denote the induced map from geodesics in $X_{i,j-1}$ to geodesics in
 $X_{ij}$. Let $\mu_{ij}\subset{X_{ij}}$
 denote the image of $\mu_i$
under the composition of the maps $\Phi_{ik}$ as $k$ runs from $1$ to
$j$. 

$\bullet$ Let $\Pi_{ij}$ denote a nearest point retraction onto
$\mu_{ij}{\subset}X_{ij}$.  On $\widetilde{E(i)} = \bigcup_jX_{ij}$, $\Pi_i$ is
  defined by 

\begin{center}
$\Pi_i{(x)} = \Pi_{ij}{(x)}$ for $x\in{X_{ij}}$.
\end{center}

$\bullet$ On $\widetilde{M_{0}} = X_{\alpha}$ (recall that $\alpha$
  corresponds to the possibly exceptional vertex of the tree of
  spaces) let $\Pi_0$ denote a nearest point retraction onto
  $B_0{({\lambda})}$.

\smallskip

On the remaining set ${\mathcal{E}}^{\prime}$, $\Pi_\lambda$ needs to
be defined with some caution. Suppose $\widetilde{E} \in
{\mathcal{E}}^{\prime}$. Observe first that no point on $\widetilde{F}
\subset \widetilde{E}$ lies at a distance less than $K$ from $B_0 (
\lambda )$, where $K$ is as chosen after Corollary
\ref{oneball2}. Also let $k$ be as in Lemma \ref{projn-bdd2}.

{\bf Case (a):} $N_{2k} ( \widetilde{G} )$ does not intersect
${\mathcal{H}}( \lambda )$. In this case, choose $y \in \widetilde{F}$
arbitrarily and define 

\begin{center}

$\Pi_E (x) = \Pi_0 (y)$ for all $x \in \widetilde{E}$.

\end{center}

For $\widetilde{E} = {\widetilde{E(i)}}$, denote $\Pi_E$ by $\Pi_i$.

{\bf Case (b):}  $N_{2k} ( \widetilde{G} )$ intersects precisely one
$\bf{H} \in {\mathcal{H}}( \lambda )$. In this case, there exists a
unique lift $\widetilde{\sigma}$ of a curve $\sigma$ on $F$,
(parabolic in $M_0$) lying on $\widetilde{F}$ at a bounded distance
from $\bf{H}$.

As in the construction of $B({\mu})$ for $\mu \subset
\widetilde{F} \subset \widetilde{E}$, construct $B (
 {\widetilde{\sigma}} ) = \bigcup_j \widetilde{\sigma}_j \subset \bigcup_j
\widetilde{S(j)}$ where  $S(j)$ denotes the $j$th member of a sequence
of truncated pleated surfaces exiting $E$ and 
$ \widetilde{\sigma}_j = \Phi_j(  \widetilde{\sigma}_{j-1} )$ ($\Phi_j$
is the map on geodesics induced by the quasi-isometry $\phi_j$ from
$\widetilde{S(j-1)}$ to  $\widetilde{S(j)}$). Note that this
construction works just as well for infinite geodesics. On
$ \widetilde{S(j)}$, define $\Pi_{\sigma j}$ to be a nearest point
projection onto  $\widetilde{\sigma}_j$. Define 

\begin{center}

$\Pi^0_\sigma (x) = \Pi_{\sigma j} (x)$ for $ x \in  \widetilde{S(j)}$.

\end{center}

Next, let $\Pi_\sigma$ denote nearest point projection of
$\widetilde{E}$ onto $ \widetilde{\sigma}$. 
Define

\begin{center}

$\Pi_E (x) = \Pi_0 \circ \Pi_\sigma (x) \circ \Pi^0_\sigma (x)$ 
for all $x \in \widetilde{E}$.

\end{center}

For $\widetilde{E} = {\widetilde{E(i)}}$, denote $\Pi_E$ by $\Pi_i$.

\smallskip

The reason for factoring the projection 
$ \Pi_\sigma (x)$ through $\Pi^0_\sigma (x)$ is that on
$\widetilde{F(i)}$ we want $\Pi_i$ to coincide with
$\Pi_\sigma$. (Else nearest point projections will have to be taken in
the $\widetilde{E(i)}$ metric and $\widetilde{F(i)}$ is not
quasiconvex in this metric, so that nearest point projections might
well differ substantially in the two metrics.)

{\bf Non-definition:} We would like to define \\
\begin{eqnarray*}
\Pi_\lambda{(x)} & = & \Pi_i{(x)},  x \in \widetilde{E(i)} \in
\mathcal{E} \\
 & = & \Pi_0{(x)},  x \in X_\alpha \\
 & = & \Pi_i{(x)},  x \in {\widetilde{E(i)}} \in
{\mathcal{E}}'
\end{eqnarray*}

\smallskip

{\bf Caveat:} Our non-definition above is not yet a real definition as
$\Pi_\lambda$ has been putatively
 defined twice on each $\widetilde{F(i)}$, once regarding  
 $\widetilde{F(i)}$ as a subset of $\widetilde{M_0}$ and again,
 regarding it as a subset of  $\widetilde{E(i)}$. We will show that
there is at most a bounded amount of discrepancy between the two
definitions and so any choice will work. So we define:

\smallskip

{\bf Definition:}\\
\begin{eqnarray*}
\Pi_\lambda{(x)} & = & \Pi_i{(x)},  x \in \widetilde{E(i)} - X_{i0} ,
  \widetilde{E(i)} \in
\mathcal{E} \\
 & = & \Pi_0{(x)},  x \in X_\alpha \\
 & = & \Pi_i{(x)},  x \in {\widetilde{E(i)}} \in
{\mathcal{E}}'
\end{eqnarray*}

\smallskip

\subsection{$\Pi_{\lambda}$ is a retract}

\smallskip

We will show in this subsection that there exists $C > 0$ such
that for all $x, y \in X = \widetilde{N_0},
d({\Pi_\lambda}(x),{\Pi_\lambda}(y)) \leq C d(x,y)$
There are three parts in the definition of $\Pi_\lambda$ as described
above. We discuss the three cases separately. 

Since we are dealing with graphs, it suffices to prove the result for
$d(x,y) = 1$. We need to check the following:

\begin{enumerate}
\item $x \in \widetilde{F(i)} =
  \widetilde{E(i)}{\cap}\widetilde{M_0}$ for some 
  $\widetilde{E(i)}\in{\mathcal{E}}$. We want to show
  that $d({\Pi_i}(x),{\Pi_0}(x)) \leq C$  \\
\item  $x, y \in \widetilde{E(i)}$ for some
  $\widetilde{E(i)}\in\mathcal{E}$ and $d(x,y)=1$. We want to show
  that $d({\Pi_\lambda}(x),{\Pi_\lambda}(y)) \leq C$ \\
\item $x, y \in \widetilde{M_0} = X_\alpha$ and $d(x,y)=1$. We want to show
  that $d({\Pi_0}(x),{\Pi_0}(y)) \leq C$ \\
\item $x \in \widetilde{F(i)} =
  \widetilde{E(i)}{\cap}\widetilde{M_0}$ for some 
  $\widetilde{E(i)}\in{\mathcal{E}}^{\prime}$. We want to show
  that $d({\Pi_i}(x),{\Pi_0}(x)) \leq C$ \\
\item  $x, y \in \widetilde{F(i)} =
  \widetilde{E(i)}{\cap}\widetilde{M_0}$ for some 
  $\widetilde{E(i)}\in{\mathcal{E}}^{\prime}$,and $d(x,y) = 1$. We want to show
  that $d({\Pi_i}(x),{\Pi_i}(y)) \leq C$. 
\end{enumerate}

In the above check-list, the first and fourth steps ensure that there is
approximate agreement on $\widetilde{F_i}$ for the two possible
definitions of $\Pi_\lambda$ occurring in the non-definition given
earlier. This ensures smooth passage from the non-definition to the
definition. The rest of the steps
 are required to prove the three cases in the
definition of $\Pi_\lambda$.

\smallskip

\noindent
{\bf Step 1:} {\underline{\bf Bounded discrepancy  on
  $\widetilde{F(i)}$ when ${\widetilde{E(i)}} \in
  \mathcal{E}$}} \\ 
 The next lemma follows
easily from the fact that local quasigeodesics in a hyperbolic metric
space
are quasigeodesics \cite{GhH} (See also Lemma 3.1 of
\cite{mitra-trees}).
  If $x, y$ are points in a hyperbolic
metric space, $[x,y]$ will denote a geodesic joining them.

\begin{lemma}
Given $\delta , C_1 > 0$,
 there exist $D, C_0$ such that  if $a, b, c, d$
are vertices of a $\delta$-hyperbolic metric space $(Z,d)$, and 
$W \subset Z$ is a $C_1$-quasiconvex set, 
with ${d}(a,W)={d}(a,b)$,
${d}(d,W)={d}(c,d)$ and ${d}(b,c)\geq{D}$
then $[a,b]\cup{W}\cup{[c,d]}$ is $C_0$-quasiconvex. Also, if
 $[b,c]_W$ denotes the shortest path in $W$ joining $b, c$, then
 $[a,b]{\cup}[b,c]_W{\cup}[c,d]$ is a $C_0$-quasigeodesic. Further, if
 $[b,c]_{amb}$ denotes an ambient quasigeodesic in $(Z - W)$
 (adding on $b, c$ as the initial and terminal points), then 
 $[a,b]{\cup}[b,c]_{amb}{\cup}[c,d]$ is an ambient $C_0$-quasigeodesic
 in $(Z-W)$.
\label{perps}
\end{lemma}

\smallskip

The next couple of Lemmas are taken from \cite{mitra-trees}, where they
  are stated in the general framework of trees of hyperbolic metric
  spaces satisfying the q.i. embedded condition.
$X_v , X_e, f_v$ denote respectively vertex space, edge space and
  q.i.  embedding of edge space in vertex space.

\begin{lemma} (Lemma 3.6 of \cite{mitra-trees}) Let $X$ be a tree(T) of hyperbolic metric spaces and $v$
  be a vertex of $T$. Let $C > 0$.
Let $\mu_1 = [a,b] \subset X_v$ be a geodesic and let $e$ be
an edge of $T$ incident on $v$. Let $p, q \in {N_C}({\mu_1}){\cap}{f_v}({X_e})$
be such that ${d_v}(p,q)$ is maximal. Let $\mu_2$ be a geodesic in $X_v$
joining $p, q$.  If $r \in {N_C}({\mu_1}){\cap}{f_v}({X_e})$, then
${d_v}(r,{\mu_2}) \leq D_1$ for some constant $D_1$ depending only
on $C,  \delta$.
\label{nbddist}
\end{lemma}

{\bf Proof:} Let $\pi$ denote a nearest point projection onto $\mu_1$.
Since $\mu_2$  and $[\pi{(p)},\pi{(q)}] \subset \mu_1$ are  geodesics
whose end-points lie at distance at most $C$ apart, there exists $C^{\prime}$
such that $[\pi{(p)},\pi{(q)}] \subset N_{C^\prime}({\mu_2})$. If
$\pi{(r)} \in [\pi{(p)},\pi{(q)}] $, then 
\begin{center}
$d(r,{\mu_2}) \leq C + C^{\prime} $.
\end{center}

If
$\pi{(r)} \notin [\pi{(p)},\pi{(q)}] $, then without loss of generality,
assume $\pi{(r)} \in [a,\pi{(p)}] \subset [a,{\pi}(q)]$. Then

\begin{eqnarray*}
d(p,q) & \geq & d(r,q)\\
         & \geq  & d({\pi}(r),{\pi}(q)) -2C \\
       &   = &   d({\pi}(r),{\pi}(p)) + d({\pi}(p),{\pi}(q))-2C \\
      &    \geq & d({\pi}(r),{\pi}(p)) + d(p,q)-4C \\
\Rightarrow  d({\pi}(r),{\pi}(p)) & \leq & 4C \\
\Rightarrow  d(r,p) & \leq & 6C \\
\Rightarrow  d(r,{\mu_2}) & \leq  & 6C.
\end{eqnarray*}

Choosing $ D_1 = max\{{C + C^{\prime}},{6C}\}$, we are through. $\Box$

\medskip

\begin{lemma} (Lemma 3.7 of \cite{mitra-trees}) Let $\mu_1$, $\mu_2$ be as in Lemma \ref{nbddist} above.
Let $\pi_i$ denote nearest point projections onto $\mu_i$ $(i= 1,2)$.
If $p\in{f_v}({X_e})$, then $d({\pi_1}(p),{\pi_2}(p)) \leq C_6$
for some constant $C_6$ depending on $\delta$ alone.
\label{proj2}
\end{lemma}

\medskip

{\bf Proof:} If $d({\pi_1}(p),{\pi_1}{\cdot}{\pi_2}(p)) \leq D$, then
$d({\pi_1}(p),{\pi_2}(p)) \leq C + D$.

Else, suppose $d({\pi_1}(p),{\pi_1}{\cdot}{\pi_2}(p)) > D$. Then
$[p, \pi_1(p)] \cup [\pi_1(p),{\pi_1}{\cdot}{\pi_2}(p)] \cup
[{\pi_1}{\cdot}{\pi_2}(p), \pi_2(p)]$ is a (uniform) quasigeodesic
by Lemma \ref{perps} . But $p, \pi_2 (p) \in {f_v}({X_e})$ which
is quasiconvex. Hence
 there exists (uniform) $C_1$ and $r \in {f_v}({X_e})$ such that
$d(r,{\pi_1}(p)) \leq C_1$, .

Then, by Lemma \ref{nbddist} above, there exists $s\in{\mu_2}$ such that
$d(s,{\pi_1}(p)) \leq C_1 + D_1$.

Again (See for instance the  proof of Lemma 3.2 of \cite{mitra-trees}), 
${(p,s)}_{{\pi_2}(p)} \leq 2\delta$.  Hence, 
${(p,{\pi_1}(p))}_{{\pi_2}(p)} \leq 2\delta + C_1 + D_1$.  

Similarly, $(p,{\pi_1}{\cdot}{\pi_2}(p))_{{\pi_1}(p)} \leq 2 \delta$.
Hence, 
${(p,{\pi_2}(p))}_{{\pi_1}(p)} \leq 2\delta + C $.  

Therefore, using the definition of Gromov inner product, \\
 $d({\pi_1}(p),{\pi_2}(p)) =
(p,{\pi_1}(p))_{{\pi_2}(p)} + {(p,{\pi_2}(p))}_{{\pi_1}(p)} 
\leq 4\delta + 2C + D_1 $.
Choosing $C_6 = 4\delta + 2C + D_1 $ we are through. $\Box$

\medskip

Though we have stated and proved Lemmas \ref{proj2} and Lemma
\ref{nbddist}  for geodesic
segments, the proof goes through for quasiconvex sets. What Lemma \ref{proj2}
effectively says is the following: \\

\smallskip

We start with a geodesic $\mu_1$ in a
 hyperbolic metric space $X =
X_v$. $W = f_v (X_e) \subset X_v$ is a $C_1$-quasiconvex set. We
consider the set of points $P = $ $\{ $ $p \in W : d(p. \mu_1) \leq
C_2$ $\}$. Choose $x, y \in P$ such that $d(x,y)$ is maximal. Let
$\mu_2$ be the geodesic joining $x, y$. Then nearest-point retractions
onto  $\mu_1, \mu_2$ almost agree for points in $W$. 

This becomes easier to grasp if all sets in sight are convex subsets
of ${\mathbb{H}}^n$. We start with $X = {\mathbb{H}}^n$. $W \subset X$
is convex. $\mu_1$ is replaced by another convex set $V \subset X$. We
consider the intersection of neighborhoods $N_\epsilon (V) \cap
N_\epsilon (W) = V_1 \neq \emptyset$. Then nearest point retractions
onto $V_1$ and $V$ almost agree for $p \in W$. This is the context of
Lemma 
\ref{proj2} in general. We quasify this below.

\begin{cor}
Given $\delta , C, C_1 > 0$, there exists $D > 0$ such that if
\begin{enumerate}
\item $A \subset Y$, $B \subset Z$, $Y \subset Z$ are inclusions of
$C_1$-quasiconvex sets into $\delta$ hyperbolic metric spaces,\\
\item $A$ and $(N_C(B)) \cap Y)$ are within a Hausdorff distance of
  $D$ of each other \\
\item $\Pi_A$ denote a nearest point projection of $Y$ onto $A$, and $\Pi_B$
denote a nearest point projection of $Z$ onto $B$, \\
\end{enumerate}
then for all $y \in
Y$, $d( \Pi_A (y), \Pi_B (y) ) \leq D$. 
\label{projcor}
\end{cor}

There exist (possibly) curves $\sigma_{ij} $ on $F(i)$  homotopic
to some curves on components $\Delta_j$ of $P$ in $(M,P)$.
Then $\sigma_{ij}$ is one boundary component of a (possibly immersed)
annulus $A_{ij}$ whose other boundary component lies on $\Delta_j$.
Let $J(i) = F(i) \cup \bigcup{A_{ij}} \subset M_0$. Lift $J(i)$ to the
universal cover $\widetilde{M_0}$ to get copies of 
$\widetilde{J(i)}$.

We will show  that on $\widetilde{F(i)}$, the two putative  definitions of
$\Pi_\lambda$  almost agree. Suppose 
$\mu_i = [{a_i},{b_i}]_F
\subset\widetilde{F(i)}\subset\widetilde{G(i)}$ is a geodesic in the
path-metric on $\widetilde{F(i)}$. 
As in the construction of $B_0({\lambda})$, we can construct a
`hyperbolic' geodesic $\mu_i^h$ joining $a_i, b_i$ in
$\widetilde{G(i)}$. Such a geodesic 
has unique entry and exit points for every horodisk and hence unique
`Euclidean' geodesics joining them on the bounding
horocycles. Replacing the hyperbolic segments by Euclidean segments,
we obtain an ambient quasigeodesic $\mu_i^a$ (By Theorem \ref{ctm}) in
$\widetilde{J(i)}$.   For $x \in \widetilde{F(i)}$, let $\Pi_{i1}$
denote a nearest point retraction onto $\mu_i$ in the path metric on
$\widetilde{F(i)}$. Also, let $\Pi_{i2}$ denote a nearest point retraction onto
${\mu_i^a}\subset\widetilde{J(i)}$ in the path metric on
$\widetilde{J(i)}$.
Then from Lemma \ref{nbddist} and Lemma \ref{proj2} above, we get

\begin{lemma}
There exists $C > 0$ such that for all $i$ and $x \in \widetilde{F(i)}
\subset \widetilde{J(i)}$, 

\begin{center}

$d({\Pi_{i1}}(x),{\Pi_{i2}}(x)) \leq C$

\end{center}
\label{closeprojns1}
\end{lemma}

\smallskip

Again, as in the construction of $B_0({\lambda})$,
each
of the Euclidean geodesic segments in $\mu_i^a$
 along with the hyperbolic geodesic
segments of $\mu_i^h$ joining its end-points (and lying entirely
within the corresponding 
horodisk) bounds 
a 
 totally geodesic disk. The union of $\mu_i^a$ and all these totally
geodesic disks is denoted by $B_0{( \mu_i)}$. There exists $C_1 >
0$ such that each $B_0{({\mu_i})}$ is $C_1$-quasiconvex in
$\widetilde{G_i}$ and hence in $\widetilde{M_0}$.
 (See for instance McMullen \cite{ctm-locconn}
Section 8.)

Let $\Pi_{i3}$ denote a nearest point projection onto  $B_0{( \mu_i)}$,
in $\widetilde{G(i)}$.

\begin{lemma}
There exists $C > 0$ such that for all $i$ and $x \in \widetilde{F(i)}
\subset \widetilde{J(i)}$, 

\begin{center}

$d({\Pi_{i2}}(x),{\Pi_{i3}}(x)) \leq C$

\end{center}
\label{closeprojns2}
\end{lemma}

{\bf Proof:} For the purposes of this Lemma,  $[a,b]$ will denote a
geodesic in $\widetilde{G(i)}$.
First, $[x,{\Pi_{i2}{(x)}}]\cup\mu_i^a$ is a tripod in $\widetilde{J(i)}$
and hence it is
uniformly quasiconvex from Lemma \ref{perps}. Again, by the same
Lemma, $[x,{\Pi_{i3}{(x)}}]\cup{B_0}({\mu_i})$ is quasiconvex. 
Further, $\Pi_{i3}{(x)} \in \mu_i^a$ (since $\mu_i^a$ forms the boundary of
${B_0}({\mu_i})$ in $\widetilde{G(i)}$ and separates it from the
$\widetilde{J(i)}$). Therefore,$[x,{\Pi_{i2}{(x)}}]\cup
[{\Pi_{i2}{(x)}},{a_i}]$ and (from Theorem \ref{ctm}) $[x,{\Pi_{i3}{(x)}}]\cup
[{\Pi_{i3}{(x)}},{a_i}]$ are both ambient quasigeodesics in the
hyperbolic metric space $\widetilde{G(i)}$. In the same way,
$[x,{\Pi_{i2}{(x)}}]\cup
[{\Pi_{i2}{(x)}},{b_i}]$ and $[x,{\Pi_{i3}{(x)}}]\cup
[{\Pi_{i3}{(x)}},{b_i}]$ are both ambient quasigeodesics in the
hyperbolic metric space $\widetilde{G(i)}$.

Hence, $\Pi_{i3}(x)$ must lie in a bounded neighborhood of 
$[x,{\Pi_{i2}{(x)}}]\cup
[{\Pi_{i2}{(x)}},{a_i}]$, as also $[x,{\Pi_{i2}{(x)}}]\cup
[{\Pi_{i2}{(x)}},{b_i}]$, as also $\mu_i^a$. Hence, it must lie close
to the intersection of these three sets, which is $\Pi_{i2}(x)$.
This proves the Lemma. $\Box$

\smallskip

Recall that $\Pi_0$ denotes nearest point projection onto
 $B_0{({\lambda})}$. Now,
 $B_0{({\mu_i})}$ lies in a bounded neighborhood of
 $B_0{({\lambda})}$, since $\mu_i$ lies in a bounded neighborhood of
 $\lambda$. 
Now from Corollary \ref{projcor}, choosing $Y = \widetilde{G(i)}$, $Z
 = \widetilde{M_0}$, $A = B_0 ( \mu_i )$,  $B = B_0 ( \lambda )$, we
 obtain

\begin{cor}
There exists $C > 0$ such that for all $x \in \widetilde{F(i)} \subset
\widetilde{J(i)}  \subset
\widetilde{G(i)}  \subset
\widetilde{M_0} $, $d ( \Pi_0 (x), \Pi_{i3} (x) ) \leq C$.
\label{projcor2}
\end{cor}

Combining Lemma \ref{closeprojns1}, Lemma \ref{closeprojns2} and
Corollary \ref{projcor2}, we get

\begin{lemma}
Recall that for $x \in \widetilde{F(i)}$, $\Pi_i$ denotes nearest point
retraction onto $\mu_i$ and  $\Pi_0$ denotes nearest point
retraction onto ${B_0}({\lambda})$. Then $d(\Pi_i{(x)},\Pi_{\alpha}(x)
\leq C_4$ for some $C_4 \geq 0$.
\label{almostagree1}
\end{lemma}

\noindent
{\bf Step 2:} {\underline{\bf Retract on  $\widetilde{E(i)}$ for 
  $\widetilde{E(i)}\in\mathcal{E}$} } \\
A number of lemmas will
be necessary. These are lifted directly from \cite{mitra-trees}. We
state them here without proof.

The following Lemma  says nearest point projections in a $\delta$-hyperbolic
metric space do not increase distances much.

\begin{lemma}
(Lemma 3.2 of \cite{mitra-trees})
Let $(Y,d)$ be a $\delta$-hyperbolic metric space
 and  let $\mu\subset{Y}$ be a $C$-quasiconvex subset, 
e.g. a geodesic segment.
Let ${\pi}:Y\rightarrow\mu$ map $y\in{Y}$ to a point on
$\mu$ nearest to $y$. Then $d{(\pi{(x)},\pi{(y)})}\leq{C_3}d{(x,y)}$ for
all $x, y\in{Y}$ where $C_3$ depends only on $\delta, C$.
\label{easyprojnlemma}
\end{lemma}

\medskip

The following Lemma says that nearest point projections and quasi-
isometries
in hyperbolic metric spaces `almost commute'.

\begin{lemma}
(Lemma 3.5 of \cite{mitra-trees}) Suppose $(Y,d)$ is $\delta$-hyperbolic.
Let $\mu_1$ be some geodesic segment in $Y$ joining $a, b$ and let $p$
be any vertex of $Y$. Also let $q$ be a vertex on $\mu_1$ such that
${d}(p,q)\leq{d}(p,x)$ for any $x\in\mu_1$. 
Let $\phi$ be a $(K,{\epsilon})$ - quasi-isometry from $Y$ to itself.
Let $\mu_2$ be a geodesic segment 
in $Y$ joining ${\phi}(a)$ to ${\phi}(b)$ . Let
$r$ be a point on $\mu_2$ such that ${d}({\phi}(p),r)\leq{d}({\phi(p)},x)$ for $x\in\mu_2$.
Then ${d}(r,{\phi}(q))\leq{C_4}$ for some constant $C_4$ depending   only on
$K, \epsilon $ and $\delta$. 
\label{cruciallemma}
\end{lemma}
\medskip

\begin{theorem}
There exists $C > 0$, such that
if  $x, y \in \widetilde{E(i)}$ for some
  $\widetilde{E(i)}\in\mathcal{E}$, and $d(x,y) = 1$, then 
 $d({\Pi_i}(x),{\Pi_i}(y)) \leq C$.
\label{ct-trees}
\end{theorem}

{\bf Proof:} 
{\bf Case (a):} {\it $x, y \in X_{ij}$ for some $i = 1, \cdots s$, $j = 0,
1, \cdots$.} \\
This follows directly from Lemma \ref{easyprojnlemma}.\\

{\bf Case (b):} $x \in X_{i,j-1}$ and $y \in X_{ij}$ 
for some $i = 1, \cdots s$, $j = 0,
1, \cdots$. \\ 
Recall that $\mu_{ik} = B(\mu_i)\cap{X_{ik}}$ for all $i, k$. Then
$\Pi_i{(z)} = \Pi_{ik}(x)$ for $z \in X_{ik}$. Also, 
from Lemma \ref{equispaced}, there exist $K, \epsilon > 0$ such that
for all $i, k$, $\phi_{ik}$ is a $(K, \epsilon)$ - quasi-isometry.
Hence, ${\phi_{ij}}({\mu_{i,j-1}})$ is a 
$({K,\epsilon})$-quasigeodesic  lying in a bounded 
 $K^{\prime}$-neighborhood 
of  ${\Phi_{ij}}({\mu_{i,j-1}}) = \mu_{ij}$
 where $K^{\prime}$ depends only on $K, {\epsilon}, \delta$. 

Now ${\phi_{ij}}\circ\Pi_{i,j-1}(x)$ lies on this quasi-geodesic.
\\By Lemma \ref{cruciallemma}, there exists $C_0 > 0$ such that
$d({\phi_{ij}}\circ\Pi_{i,j-1}(x),\Pi_{ij}\circ{\phi_{i,j-1}}(x)) \leq
C_0$.
\\Also $d(x,y) = 1 = d(x, \phi_{ij}(x))$. Hence,
 $d(\phi_{ij}(x), y) \leq 2$ and $\phi_{ij}(x), y \in
X_{ij}$. Therefore, by \ref{easyprojnlemma},  there exists $C_1$ such
that
$d(\Pi_{ij}\circ\phi_{ij}(x), \Pi_{ij}(y)) \leq 2C_1$.
\\Hence
\begin{eqnarray*}
{d}(\Pi_i(x),\Pi_i(y)) & = & d(\Pi_{i,j-1}(x),\Pi_{ij}(y)) \\
 & \leq & d(\Pi_{i,j-1}(x),{\phi_{ij}}\circ\Pi_{i,j-1}(x)) + \\
  &      & d({\phi_{ij}}\circ\Pi_{i,j-1}(x),\Pi_{ij}\circ{\phi_{i,j-1}}(x)) +
d(\Pi_{ij}\circ\phi_{ij}(x), \Pi_{ij}(y))\\
 & \leq & 1 + C_0 + 2C_1
\end{eqnarray*}

Choosing $C = (1 + C_0 + 2C_1)$, we are through. $\Box$

\smallskip

\noindent 
{\bf Step 3:} {\underline{\bf Retract on $ \widetilde{M_0} = X_\alpha$}} \\
This case follows directly from Lemma \ref{easyprojnlemma}. We state
 a special case required for our specific purposes.

\begin{lemma}
There exists $C > 0$ such that if 
 $x, y \in \widetilde{M_0} = X_\alpha$ and $d(x,y) = 1$, then
 $d({\Pi_0}(x),{\Pi_0}(y)) \leq C$.
\label{projn-rootvertex}
\end{lemma}

\noindent  {\bf Step 4:} {\underline {\bf Bounded discrepancy on
  $\widetilde{F(i)}$ when ${\widetilde{E(i)}} \in
  {\mathcal{E}}^{\prime}$}} \\
We want to show that for 
 $x, y \in \widetilde{F(i)} \subset {\widetilde{E(i)}} \in
  {\mathcal{E}}^{\prime}$, $d(\Pi_0 (x), \Pi_i (x))$ is uniformly
 bounded. 

Recall from the definition of $\widetilde{G(i)}$, two case may arise:
\\ 
{\bf Case (a):} $N_{2k} ( \widetilde{G(i)} )$ does not intersect
${\mathcal{H}}( \lambda )$. In this case, recall that $y \in \widetilde{F(i)}$
is chosen arbitrarily, and we had defined

\begin{center}

$\Pi_i (x) = \Pi_0 (y)$ for all $x \in \widetilde{E}$.

\end{center}

Next, by the choice of 
$K = 3k + 4\delta + k^\prime$ 
in the definition of ${\mathcal{E}}^{\prime}$, and Corollary
\ref{oneball2}, we have that the diameter of the set $\Pi_0 (
\widetilde{G(i)} )$ and hence  $\Pi_0 (
\widetilde{F(i)} )$ is bounded by $k$ (from Lemma \ref{projn-bdd2} ).

Hence, $d(\Pi_0 (x), \Pi_i (x)) = d(\Pi_0 (x), \Pi_0 (y)) \leq k$.

\noindent 
{\bf Case (b):}  $N_{2k} ( \widetilde{G(i)} )$ intersects precisely one
$\bf{H} \in {\mathcal{H}}( \lambda )$. In this case, there exists a
unique lift $\widetilde{\sigma}$ of a curve $\sigma$ on $F(i)$,
(parabolic in $M_0$) lying on $\widetilde{F(i)}$ at a bounded distance
from $\bf{H}$. Recall that  $\Pi_\sigma$ denotes nearest point projection of
$\widetilde{E(i)}$ onto $ \widetilde{\sigma}$ and $\Pi_\sigma^0 (x)$
denotes nearest point projection onto $B( \widetilde{\sigma })$.
On $\widetilde{F(i)}$, $\Pi_\sigma^0$ coincides with $\Pi_{\sigma} \circ
\Pi^0_\sigma$. Also recall that

\begin{center}

$\Pi_i (x) = \Pi_0 \circ \Pi_\sigma  \circ \Pi^0_\sigma (x) $ for all 
  $x \in \widetilde{E(i)}$. 

\end{center}

On $\widetilde{F(i)}$, $\Pi_\sigma^0$ coincides with $\Pi_\sigma \circ
\Pi^0_\sigma$.
Thus, we want to show that 
$d(\Pi_0 (x),\Pi_0 \circ \Pi^0_\sigma (x))$ is uniformly bounded on
$\widetilde{F(i)}$.

Here,  $N_k ( \widetilde{G(i)} ) \cap N_k (B_0
(\lambda ))$ lies in some horoball $\bf{H} \in \mathcal{H} ( \lambda
)$ while $N_k ( \widetilde{F(i)} ) \cap N_k (B_0
(\lambda )) = \emptyset$ (as per definition of ${\mathcal{E}}^{\prime}$).
By Lemma \ref{projn-bdd2} $\Pi_0 ( \widetilde{F(i)} )  \subset N_k (
\widetilde{(A \cup \kappa )}\subset N_k (
\bf{H} )$, where $\kappa \subset G(i)$ is the cusp whose boundary
$\eta$ bounds the annulus $A$ along with $\sigma$ in $M_0$. 
(Recall from the discussion following Lemma \ref{projn-bdd2}
that each $\eta$ forms the boundary of
 a totally geodesic 2 dimensional subset 
$\kappa$ of a $Z$ cusp or a $Z+Z$-cusp of $M_0$. ) Also, 
$\widetilde{(A \cup \kappa )} = \widetilde{A_\kappa } $ (say) is
$C_3$-quasiconvex for some $C_3$.

 Let $\Pi_\kappa$ denote nearest
point projection of $\widetilde{G(i)}$ onto
$\widetilde{A_\kappa}$ in the metric on $\widetilde{M_0}$. Also, let
$\Pi_0^{aug}$ denote nearest point 
projection onto $B_0^{aug} ( \lambda )$. Then by Corollary
\ref{projcor}
there exists $C_4 > 0$ such that 

\begin{center}

$d( \Pi_\kappa (x), \Pi_0^{aug} (x) )
\leq C_4$ for all $x \in \widetilde{G(i)}$.

\end{center}

Next, $B_0 ( \lambda ) \subset B_0^{aug} ( \lambda ) \subset
\widetilde{M_0}$ and both $B_0 ( \lambda )$ and $B_0^{aug} ( \lambda
)$ are $C_1$-quasiconvex. Both $\Pi_0$ and $\Pi_0 \circ \Pi_0^{aug}$
are large-scale Lipschitz retracts onto $B_0 ( \lambda )$, the only difference
being that the latter factors through a large-scale Lipschitz retract
$\Pi_0^{aug}$ onto $B_0^{aug} ( \lambda
)$. Hence the two maps must send $x$ to close by points, i.e.
there exists $C_5 > 0$ such that 

\begin{center}

 $d( \Pi_0 \circ \Pi_0^{aug} (x), \Pi_0 (x) )
\leq C_5$ for all $x \in \widetilde{M_0}$.

\end{center}

Using the above two inequalities along with Lemma \ref{easyprojnlemma}
we get
a $C_6 > 0$ such that

\begin{center}

 $d( \Pi_0 \circ \Pi_\kappa (x), \Pi_0 (x) )
\leq C_6$ for all $x \in \widetilde{M_0}$.

\end{center}

Finally, we observe that for $x \in \widetilde{F(i)}$, since $\Pi^0_\sigma$
denotes nearest point projection onto $\widetilde{\sigma }$, then 
$d( \Pi_\kappa (x), \Pi^0_\sigma (x) )$ is bounded. This is because
$\widetilde{\sigma }$ separates $\widetilde{G(i) }$ into
$\widetilde{F(i) }$
and $\widetilde{A_\kappa }$. Hence, the geodesic joining $x$ to
$\Pi_\kappa (x)$ in $\widetilde{G(i) }$ must cut $\widetilde{\sigma }$
at some point which therefore must coincide with 
$\Pi^0_\sigma (x)$, provided we take nearest point projections in
$\widetilde{G(i)}$. 
Since $\widetilde{G(i)}$ is quasiconvex in $\widetilde{M_0}$,
therefore by Corollary \ref{projcor}, we might as well take
projections in the $\widetilde{M_0}$ metric and we have a $C_7 > 0$ such that

\begin{center}

 $d( \Pi_\kappa (x), \Pi^0_\sigma (x) )
\leq C_7$ for all $x \in \widetilde{F(i)}$.

\end{center}

Combining the above two inequalities and using Lemma
\ref{easyprojnlemma} for $\Pi_0$, we obtain finally,

\begin{lemma}
There exists
$C > 0$ such that

\begin{center}

 $d( \Pi_0 \circ \Pi^0_\sigma (x), \Pi_0 (x) )
\leq C$ for all $x \in \widetilde{F(i)}$.

\end{center}

Since on such a $\widetilde{F(i)}$, $\Pi_0 \circ \Pi^0_\sigma$ is
denoted by $\Pi_i$, we have 
 $d( \Pi_i (x), \Pi_0 (x) )
\leq C$ for all $x \in \widetilde{F(i)}$.
\label{almostagree2}
\end{lemma}

\smallskip

The three results Lemma \ref{almostagree1}, Case (a) of Step 4 above
and Lemma \ref{almostagree2} show that $\Pi_0$ and $\Pi_i$ almost
agree, i.e. have at most a bounded amount of discrepancy for all
$\widetilde{F(i)}$. We summarize these three in the following useful
proposition.

\begin{prop}
There exists $C > 0$ such that for all $\widetilde{F(i)}$ and all 
 $x \in \widetilde{F(i)}$,

\begin{center}

 $d( \Pi_i (x), \Pi_0 (x) ) \leq C$.

\end{center}
\label{almostagree}
\end{prop}

\noindent {\bf Step 5:} {\underline {\bf Retract on  $\widetilde{E(i)}$ for 
  $\widetilde{E(i)}\in{\mathcal{E}}^{\prime}$ }} \\
$\Pi_i (x) = \Pi_0 \circ \Pi_\sigma  \circ \Pi^0_\sigma (x) $ for all 
  $x \in \widetilde{E(i)}$. 
By the proof of 
Theorem \ref{ct-trees}, there exists $C_1 > 0$ such that for all
  $x, y \in \widetilde{E(i)} \in {\mathcal{E}}^{\prime}$ with $d(x,y) = 1$, we have

\begin{center}

$d( \Pi^0_\sigma (x) ,  \Pi^0_\sigma (y)) \leq C_1$.

\end{center}

(Though Theorem \ref{ct-trees} is stated for 
$\widetilde{E(i)}\in{\mathcal{E}}$, all that we need
for the above assertion is the intrinsic metric on
$\widetilde{E(i)}$ and here the proof is the same.)

Also, by Lemma \ref{easyprojnlemma}, there exists $C_2 > 0$ such that
for $x, y \in \widetilde{M_0}$,

\begin{center}

$d( \Pi_0 (x) ,  \Pi_0 (y)) \leq C_2 d(x,y)$.

\end{center}

Hence, it would be enough to prove that there exists $C_3 > 0$ such
that

\begin{center}

$d( \Pi_\sigma (x) ,  \Pi_\sigma (y)) \leq C_3 d(x,y)$

\end{center}
for all $x, y \in \widetilde{E(i)}$.

For this, it would suffice to show that $\widetilde{\sigma}$ is  quasiconvex in
$ \widetilde{E(i)}$. Recall that there exists a simply degenerate
end  $E^h$, one
of whose lifts to the universal cover is 
$ \widetilde{E^h(i)}$. $ \widetilde{E^h(i)}$ with cusps removed is
$ \widetilde{E(i)}$. We shall show that $\widetilde{\sigma}$ is
quasiconvex in $ \widetilde{E^h(i)}$. We give $E^h$ a convex
hyperbolic structure by taking the cover of $N^h$ corresponding
to $\pi_1(E^h)$ and looking at its convex core.
This makes the universal cover of 
$E^h$ (equipped with such a convex
hyperbolic structure) quasi-isometric
to
$ \widetilde{E^h(i)}$. For the purposes of this Step, therefore, we
assume
that $ \widetilde{E^h(i)}$ is the universal cover of a simply
degenerate hyperbolic manifold (so as to avoid the extra expository
complication of mapping via quasi-isometries). Now, $\sigma \subset
E^h$ is a closed curve, and is therefore homotopic by a bounded
homotopy to a closed geodesic $\sigma^{\prime}$. Then any lift 
 $\widetilde{\sigma^{\prime}}$ is a hyperbolic geodesic in 
the convex hyperbolic manifold  $ \widetilde{E^h(i)}$. 
Since  $\widetilde{\sigma}$ lies at a bounded distance from
$\widetilde{\sigma^{\prime}}$, we conclude that  $\widetilde{\sigma}$
is quasiconvex in $ \widetilde{E^h(i)}$ and hence there exists a
nearest point projection of $ \widetilde{E^h(i)}$ onto
$\widetilde{\sigma}$
which stretches distances by a bounded factor by Lemma
\ref{easyprojnlemma}. 
Restricting this map to $ \widetilde{E(i)}$, the retraction property
persists. Since the path metric on $ \widetilde{E(i)}$ dominates the
metric
induced from $ \widetilde{E^h(i)}$, and since the metric on $
\widetilde{\sigma}$ remains undisturbed on removing cusps, 
we find  (reverting to our description of spaces as graphs)
that there exists $C_3 > 0$ such
that

\begin{center}

$d( \Pi_\sigma (x) ,  \Pi_\sigma (y)) \leq C_3 d(x,y)$

\end{center}
for all $x, y \in \widetilde{E(i)}$.

Combining the three above equations, we conclude:

\begin{prop}
There exists $C > 0$, such that
if  $x, y \in \widetilde{E(i)}$ for some
  $\widetilde{E(i)}\in {\mathcal{E}}^{\prime}$,
 and $d(x,y) = 1$, then 
 $d({\Pi_i}(x),{\Pi_i}(y)) \leq C$.
\label{retn-3rdcase}
\end{prop}

\noindent {\underline {\bf $\Pi_\lambda$ is a retract: Proof}} \\

\smallskip

We combine the five steps above to conclude that $\Pi_\lambda$ is a
retract. Recall the definition of $\Pi_\lambda$.\\

\smallskip

{\bf Definition:}\\
\begin{eqnarray*}
\Pi_\lambda{(x)} & = & \Pi_i{(x)},  x \in \widetilde{E(i)} - X_{i0} ,
  \widetilde{E(i)} \in
\mathcal{E} \\
 & = & \Pi_0{(x)},  x \in X_\alpha \\
 & = & \Pi_i{(x)},  x \in \widetilde{E(i)} - X_{i0} ,{\widetilde{E(i)}} \in
{\mathcal{E}}'
\end{eqnarray*}

Also recall 
  that $\widetilde{N_0}$ is quasi-isometric to a graph $X$
which can be regarded as a tree of hyperbolic metric spaces with
possibly exceptional vertex $\alpha$ satisfying the q.i. embedded
condition. (See Lemma \ref{equispaced}, Lemma \ref{hyptree}
 and the discussion
 preceding Lemma \ref{equispaced} where this is elaborated.)

We are now in a position to prove our main technical theorem.

\begin{theorem}
There exists $C_0\geq{0}$ such that 
${d}({\Pi_\lambda}(x),{\Pi_\lambda}(y))\leq C_0{d}(x,y)$ for $x, y$
vertices of $X$.
\label{mainref}
\end{theorem}

{\it Proof:} It suffices to prove 
the theorem when ${d}(x,y)=1$.

\begin{enumerate}
\item For $x, y \in \widetilde{E(i)} - \widetilde{M_0}$, 
$ \widetilde{E(i)} \in \mathcal{E}$, this follows from Theorem
  \ref{ct-trees} in Step 2 above. \\
\item For $x, y \in X_\alpha = \widetilde{M_0}$, this follows
  from Lemma \ref{projn-rootvertex} in Step 3 above.\\
\item  For $x, y \in \widetilde{E(i)} - \widetilde{M_0}$, 
$ \widetilde{E(i)} \in {\mathcal{E}}^{\prime}$, this follows from Proposition
  \ref{retn-3rdcase} in Step 5 above. \\
\item $x \in X_\alpha , y \in X_{i1} \subset \widetilde{E(i)}$ and
  $d(x,y) = 1$ for $\widetilde{E(i)} \in \mathcal{E} \cup
  {\mathcal{E}}^{\prime}$.  \\
\end{enumerate}
Here, 
\begin{eqnarray*}
d( \Pi_\lambda (x), \Pi_\lambda (y) ) & = & d( \Pi_0 (x), \Pi_i (y) )
\\
 & \leq & d( \Pi_0 (x), \Pi_i (x) ) + d( \Pi_i (x), \Pi_i (y) )
\end{eqnarray*}

Choose constants $C_1$, $C_2$ and $C_3$ from Proposition
\ref{almostagree}, Theorem \ref{ct-trees} and Proposition
\ref{retn-3rdcase}. Let $C = C_1 + max( C_2, C_3 )$. 
Then we get

\begin{center}

$d( \Pi_\lambda (x), \Pi_\lambda (y) ) \leq C$.
\end{center}
This proves the result.
$\Box$

\section{Proof of Main Theorem}

Recall that we started off with a geodesic $\lambda^h$ on
$\widetilde{M_{gf}}$ outside an $n$-ball around a fixed reference
point $p$, joining $a,b$. Replacing the maximal segments of
$\lambda^h$ lying inside 
$Z$-horoballs by `Euclidean' geodesics lying on the corresponding
horosphere, we obtained an ambient quasigeodesic $\lambda$ in
$\widetilde{M_0}$  joining $a, b$. $\lambda$ agrees with $\lambda^h$ off
$Z$-horoballs.
Let $\beta^h$ be the geodesic in $\widetilde{N^h}$ joining $a,
b$. Performing the same operation on $\beta^h$ for $Z$-horoballs in
$\widetilde{N^h}$, we obtain an ambient geodesic 
 $\beta_{amb}^0$ 
 in $\widetilde{N_{0}}$ joining the end-points of
$\lambda$.
We have proved the existence of a retraction $\Pi_\lambda$ in the
preceding section (Theorem \ref{mainref}. Project  $\beta_{amb}^0$ 
  onto $B_\lambda$, using $\Pi_\lambda$ to get another ambient quasigeodesic
$\beta_{amb}$. Thus, $\beta_{amb} = \Pi_\lambda ( \beta_{amb}^0
  )$.

Our starting point for this section is therefore the ambient
quasigeodesic
$\beta_{amb} \subset B_\lambda \subset \widetilde{N_0}$.

\subsection{Quasigeodesic rays}

\medskip

The purpose of this subsection is to construct for any $\widetilde{E(i)}
\in \mathcal{E}$ and $x \in B( \mu_i ) = B_\lambda \cap
\widetilde{E(i)}$ a quasi-geodesic ray $r_x \subset B_\lambda$ passing
through $x$. $r_x$ can be regarded both as a function $r_x : \{ 0 \}
\cup {\Bbb{N}} \rightarrow B( \mu_i )$ or as a subset of $ B( \mu_i )$
(the image of the function $r_x$). If $x$ lies away from cusps, so
will $r_x$. 

Fix  $\widetilde{E}
\in \mathcal{E}$. Recall that  \\
$\bullet$ $\widetilde{E} \cap \widetilde{M_0} =
\widetilde{F}$, \\
$\bullet$  $\widetilde{F} \cap B_\lambda = \mu =\mu_0$, \\
$\bullet$ $\widetilde{E} = \bigcup_i  \widetilde{S(i)} $, \\
$\bullet$ $B (\mu ) = \cup_i \mu_i$ \\
$\bullet$ $\mu_i = \Phi_i ( \mu_{i-1} ) \subset \widetilde{S(i)}$ \\

\smallskip

Also note that $\Phi_i$ is the map on geodesics induced by the
quasi-isometry $\phi_i : \widetilde{S(i-1)} \rightarrow
\widetilde{S(i)}$. The quasi-isometry $\phi_i$ is induced by a
cusp-preserving homeomorphism of truncated pleated surfaces. So we can
assume that $\phi_i$ restricted to each horocycle boundary component
of $\widetilde{S(i-1)}$ is an orientation-preserving isometry. Since
each such horocycle may naturally be identified with $\Bbb{R}$, we
might as well assume that the map on horocycles induced by $\phi_i$ is
the identity map. Now, extend $\phi_i$ to a map $\phi_i^h :
\widetilde{S^h(i-1)} \rightarrow \widetilde{S^h(i)}$ by demanding that
$\phi_i^h$ restricted to each horodisk is an isometry. Note that
$\phi_i^h$ is thus an equivariant quasi-isometry from ${\Bbb{H}}^2$ to
itself, taking horodisks isometrically to horodisks. Next, let
$\mu_i^h$ be the hyperbolic geodesic joining the end-points of
$\mu_i$. As usual, let $\mu_i^a$ denote the ambient quasigeodesic
obtained by replacing hyperbolic segments in horodisks by `Euclidean'
segments on the boundary horocycles. Then $\mu_i^a$ lies in a uniform
(independent of $i$) $C_0$-neighborhood of $\mu_i$ as they are both ambient
quasigeodesics in $\widetilde{S(i)}$. Let 

\begin{center}

$B^a ( \mu ) = \cup_i \mu_i^a $

\end{center}

Then $B^a ( \mu )$ lies in a $C_0$-neighborhood of $B( \mu )$. Let
$\mu_i^c$ denote the union of the segments of $\mu_i^a$ which lie
along horocycles and let $\mu_i^b = \mu_i^a - \mu_i^c$. Let

\begin{center}

$B^c ( \mu ) = \cup_i \mu_i^c $ \\
$B^b ( \mu ) = \cup_i \mu_i^b $\\

\end{center}

We want  to show that for all $x \in B^b ( \mu )$ there exists a
$C$-quasigeodesic $r_x : \{ 0 \} \cup \Bbb{N} \rightarrow B^b ( \mu )$
such that $x \in r_x ( \{ 0 \} \cup \Bbb{N} )$ and $r_x (i) \in
\mu_i^b$. Suppose $x \in \mu^b_k \subset B^b ( \mu )$. We define $r_x$
by starting with $r_x (k) = x$ and construct $r_x (k-i)$ and $r_x
(k+i)$ inductively (of course $(k-i)$ stops at $0$ but $(k+i)$ goes on
to infinity). For the sake of concreteness, we prove the existence of
such a $r_x (k+1)$. The same argument applies to $(k-1)$ and
inductively for the rest.

\begin{lemma}
There exists $C > 0$ such that if $r_x (k) = x \in \mu_k^b$ then there
exists $x^{\prime} \in \mu^b_{k+1}$ such that $d(x, x^{\prime} ) \leq
C$. We denote $r_x (k+1) = x^{\prime}$.
\label{qgray1}
\end{lemma}

{\bf Proof:} Let $[a,b]$ be the maximal connected component of
$\mu_k^b$  on which $x$ lies. Then there exist two horospheres
$\bf{H}_1$ and $\bf{H}_2$ such that  $a \in \bf{H}_1$ (or is the
initial point of $\mu_k$ ) and  $b \in \bf{H}_2$ (or is the terminal
point of $\mu_k$ ). Note that $[a,b]$ does not intersect any of the
horodisks of $\widetilde{S^h_k}$. Since $\phi^h_{k+1}$ preserves
horodisks, $\phi_{k+1} (a)$ lies on a horocycle (or is the initial
point of $\mu_{k+1}$ ) as does $\phi_{k+1} (b)$ (or is the terminal
point of $\mu_{k+1}$). Further, the image of $[a,b]$ is a hyperbolic
quasigeodesic (which we now denote as $\phi_{k+1} ([a,b])$) lying
outside horoballs. Let $\Phi^h_{k+1} ([a,b])$ denote the hyperbolic
geodesic joining $\phi_{k+1} (a)$ and $\phi_{k+1} (b)$.
Let $\Phi_{k+1} ([a,b])$ denote the ambient 
geodesic in $\widetilde{S(k+1)}$ 
joining $\phi_{k+1} (a)$ and $\phi_{k+1} (b)$. Therefore
$\Phi^h_{k+1} ([a,b])$ lies in a bounded neighborhood of  $\phi_{k+1}
([a,b])$ (which in turn lies at a bounded distance from
$\Phi_{k+1}([a,b])$) and hence by Theorem \ref{farb} there exists an
upper bound on 
how much $\Phi^h_{k+1} ([a,b])$ can penetrate horoballs, i.e. there
exists $C_1 > 0$  such that for all $z \in \Phi^h_{k+1} ([a,b])$,
there exists $z^{\prime} \in \Phi^h_{k+1} ([a,b])$ lying outside
horoballs with $d(z, z^{\prime}) \leq C_1$. Further, since
$\phi_{k+1}$ is a quasi-isometry     there exists $C_2 > 0$ such that
$d( \phi_{k+1} (x), \Phi^h_{k+1} ([a,b])) \leq C_2$. Hence there
exists $x^{\prime} \in  \Phi^h_{k+1} ([a,b])$ such that \\
$\noindent \bullet$ $d( \phi_{k+1} (x) , x^{\prime}) \leq C_1 + C_2$
\\
$\noindent \bullet$ $x^{\prime}$ lies outside horoballs. \\

Again,   $\Phi_{k+1} ([a,b])$
 lies at a uniformly bounded distance $\leq C_3$ from $\mu_{k+1}$ and
 so, if $c, d \in \mu_{k+1}$ such that $d(a,c) \leq C_3$ and $d(b,d)
 \leq C_3$ then the segment $[c,d]$ can penetrate only a bounded
 distance into any horoball. Hence there exists $C_4 > 0$ and
$x^{\prime \prime} \in
 [c,d]$ such that \\
$\noindent \bullet$ $ d( x^{\prime}, x^{\prime \prime} ) \leq C_4$ \\
$\noindent \bullet$ $ x^{\prime \prime}$ lies outside horoballs. \\

Hence, $d(\phi_{k+1} (x),  x^{\prime \prime}) \leq C_1 + C_2 +
C_4$. Since $d(x, \phi_{k+1} (x)) = 1$, we have, by choosing
$r_{k+1}(x) =  x^{\prime \prime}$,

\begin{center}

$d( r_k(x), r_{k+1}(x)) \leq 1 + C_1 + C_2 + C_4$.

\end{center}

Choosing $C =  1 + C_1 + C_2 + C_4$, we are through. $\Box$

\medskip

Using Lemma \ref{qgray1} repeatedly (inductively replacing $x$ with
$r_x (k + i)$ we obtain the values of $r_x (i)$ for $i \geq k$. By an
exactly similar symmetric argument, we get $r_x (k-1)$ and proceed
down to $r_x (0)$. Now for any $i$, $z \in \widetilde{S(i)}$ and
 $y \in \widetilde{S(i+1)}$, $d(z, y) \geq 1$. Hence, for any 
 $z \in \widetilde{S(i)}$ and
 $y \in \widetilde{S(j)}$, $d(z, y) \geq |i-j|$. This gives

\begin{cor} 
There exist $K, \epsilon > 0$ such that for all $x \in \mu_k^b \subset
B^b ( \mu )$ there exists a $(K, \epsilon )$ quasigeodesic ray $r_x$
such that $r_x (k) = x$ and $r_x (i) \in \mu_i^b$ for all $i = 0, 1,
2, \cdots$.
\label{qgray2}
\end{cor}

\subsection{ p-incompressible boundary}

\medskip

Till this point we have not used the hypothesis of {\it
  p-incompressibility}. We need a modification of Corollary
  \ref{qgray2} above to go down from a point in $\mu^b = \mu_0^b$ 
  to a point in $\lambda^b$.

We can change $B_\lambda$ to $B_\lambda^{a}$ by replacing each $B
( \mu ) \subset \widetilde{E(i)}$ by $B^a ( \mu )$ for $\widetilde{E(i)} \in
\mathcal{E}$. Recall that we had started off with $\lambda$ being an
ambient quasigeodesic in $\widetilde{M_0}$ constructed from
$\lambda^h$ by replacing hyperbolic geodesic segments within
$Z$-horoballs by `Euclidean' geodesic segments along horospheres. 
So, there is no need to freshly construct a $\lambda^a$ (we can think that
$\lambda = \lambda^a$). But we do need to construct $\lambda_{cc}$ by
performing the same modifications for $(Z+Z)$ - horoballs.
Just as in the construction of $\mu^a_i ,
\mu^c_i, \mu^b_i$, we can define 
$\lambda^c$ to be the collection of subsegments of $\lambda_{cc}$
lying along horoballs ( $Z$ or $(Z+Z)$), and $\lambda^b = \lambda_{cc} -
\lambda^c$. Adjoining $\lambda_{cc}$, $  
\lambda^c$ and $\lambda^b$ to $B^a ( \mu )$, $B^c (\mu)$ and $B^b (
\mu )$  we get respectively,  
$B^a (\mu , \lambda )$, $B^c (
\mu , \lambda )$ and  $B^b ( \mu , \lambda)$. Then we can extend
Corollary \ref{qgray2} above so that $\lambda$ is included.

\smallskip

For the convenience of the reader, we summarize the notation for the
different types of geodesics and quasigeodesics that we have
introduced and will use henceforth: \\
$\bullet$ $\lambda^h$ $= $ hyperbolic geodesic in $\widetilde{M_{gf}}$
joining $a, b$
\\
$\bullet$ $\lambda = \lambda^a $ $ = $ ambient quasigeodesic in
$\widetilde{M_0}$ constructed from $\lambda^h \subset
\widetilde{M_{gf}}$ \\ 
$\bullet$ $\lambda_{cc}$  $ = $ ambient quasigeodesic in
$\widetilde{M_{cc}}$ constructed from $\lambda \subset
\widetilde{M_{0}}$ \\
$\bullet$ $\lambda^c$ $ = $ part of $\lambda_{cc}$ lying along
$Z$-cusps or $Z+Z$-cusps. \\
$\bullet$ $\lambda^b$ $ = $ $\lambda_{cc} - \lambda^c$ \\
$\bullet$ $\mu_i$ $ = $ ambient geodesic in
$\widetilde{S(i)}$  \\
$\bullet$ $\mu_i^a$ $ = $  ambient quasigeodesic in
$\widetilde{S(i)}$  constructed from $\mu_i$ \\
$\bullet$ $\mu_i^c$ $ = $ part of $\mu_i^a$ lying along $Z$-cusps. \\
$\bullet$ $\mu_i^b$ $ = $ $\mu_i^a - \mu_i^c$ \\

\smallskip

 Recall that in the construction of $B_1 ( \lambda
)$ from $B_0 
( \lambda )$ we  construct $\mu$ from $\lambda$ by taking a $K$ (as in
the discussion following Corollary \ref{oneball2} ) and choosing $p, q
\in \widetilde{F} \cap N_K ( B_0 ( \lambda ) )$ with $d(p,q)$ maximal.

\smallskip

\noindent {\bf  Now suppose $(M,P)$ is  p-incompressible.} (This
assumption has effect from now till the end of Section 5.4, unless
explicitly stated otherwise.)

\smallskip

\noindent Then $\widetilde{F}$ is
quasiconvex in the hyperbolic metric space  $\widetilde{M_0}$. Hence,
$\mu$ is a quasigeodesic in 
$\widetilde{M_0}$ and therefore lies in a $K_1$-neighborhood of
$B_0 (\lambda )$ (where $K_1$ depends on $K$). Further, since $\mu$
does not penetrate horoballs, the hyperbolic geodesic $\mu^{\prime}$ in
$\widetilde{M_0}$ joining the end-points of $\mu$ can penetrate
horoballs only for some bounded length $D_1$. Also, 
there is a subsegment $\lambda_\mu = [p',q'] \subset \lambda$ such that
$d(p, p') \leq K_1$ and $d(q, q') \leq K_1$. By Theorem \ref{farb},
$\mu^{\prime}$ and $\lambda_\mu$ must have similar intersections with
horoballs, i.e.
there exists $C_0$ such that \\
$\bullet$ If only one of 
$\mu^{\prime}$
 and $\lambda_\mu$ penetrates or travels along the boundary of a horoball
  $\bf{H}$, then it can do so for a distance $ \leq C_0$. \\
$\bullet$ If both $\mu^{\prime}$ and $\lambda_\mu$ enter (or leave) a horoball
  $\bf{H}$ then their entry (or exit) points are at a distance of at
  most $C_0$ from each other.\\

Then as in Corollary \ref{qgray2}
if  $x \in \mu^b$.
 there exists $C_1 > 0$, $x^{{\prime}{\prime}} \in
 \mu^{\prime}$ and
$x^{\prime} \in
 \lambda^b$ such that 

\begin{center}

 $ d( x, x^{{\prime}{\prime}}) \leq C_1$ \\
 $ d( x^{\prime}, x^{{\prime}{\prime}}) \leq C_1$ \\

\end{center}

Hence,   $ d( x^{\prime}, x) \leq 2C_1$ \\
We have thus shown, (using Corollary \ref{qgray2} )

\begin{cor}
Suppose $(M,P)$ is a pared manifold with {\it p-incompressible boundary.} 
There exist $K, \epsilon > 0$ such that for all $x \in \mu_k^b \subset
B^b ( \mu , \lambda )$ there exists a $(K, \epsilon )$ quasigeodesic
ray $r_x : \{ -1, 0 \} \cup \Bbb{N} \rightarrow B^b ( \mu , \lambda )$
such that $r_x (k) = x$, $r_x (i) \in \mu_i^b$  for all $i = 0, 1,
2, \cdots$ and 
 $r_x (-1) \in \lambda^b$. 

\label{qgray3}
\end{cor}

\smallskip

{\bf Note:} In the discussion preceding Corollary \ref{qgray3} above,
all that we really required was that $\mu$ be a quasigeodesic in
$\widetilde{M_0}$. Quasiconvexity of $\widetilde{F}$ (following from
p-incompressibility) ensured this. We therefore state the more
general version below, as we shall require it to prove our main
Theorem \ref{main2} where p-incompressibility is relaxed to
incompressibility.

\begin{cor}
Suppose $(M,P)$ is a pared manifold with {\it incompressible
  boundary.} Given $D, \delta$ there exist $K, \epsilon$ such that the
  following holds: \\
Suppose $\mu$ is a $(D, \delta )$ hyperbolic quasigeodesic in
  $\widetilde M_0$ lying on $\widetilde F$ for some $F$. 
Then for all  $x \in \mu_k^b \subset
B^b ( \mu , \lambda )$ there exists a $(K, \epsilon )$ quasigeodesic
ray $r_x : \{ -1, 0 \} \cup \Bbb{N} \rightarrow B^b ( \mu , \lambda )$
such that $r_x (k) = x$, $r_x (i) \in \mu_i^b$  for all $i = 0, 1,
2, \cdots$ and 
 $r_x (-1) \in \lambda^b$. 

\label{qgray4}
\end{cor}

The hypothesis of Corollary \ref{qgray4} is satisfied if $\mu$ moves
along boundaries of horoballs for uniformly bounded distances. 
Equivalently, this
is satisfied if the hyperbolic geodesic $\mu^h$ joining the end-points of
$\mu$ penetrates horoballs by a uniformly bounded amount. Actually,
we do not need even this much. $\mu^h$ should penetrate {\it
  exceptional
horoballs} (See definition below) by a bounded amount. \\

\smallskip

{\bf Definition:} A cusp of $M_{gf}$ is said to be {\bf exceptional}
if there exist closed curves carried by the cusp (i.e. lying on its
boundary horocycle or horosphere) which are homotopic to
non-peripheral curves on some other boundary component of $M_{gf}$.\\
{\bf Exceptional horoballs} are lifts of exceptional cusps. \\
A geodesic is said to {\bf penetrate a  horoball $H$  by at most $D$} if any
subsegment of it lying inside {\bf $H$} has length less than or equal
to $D$.

\smallskip

The sufficient condition for the hypothesis of  Corollary
\ref{qgray4} to hold is set forth in the following Lemma.

\begin{lemma}
Let $F$ be some boundary component of $M_{0}$ and $\mu \subset
\widetilde{F}$ be a geodesic segment in the path metric on
$\widetilde{F}$. Let $\mu^h$ be the  hyperbolic geodesic in
$\widetilde{M_{gf}}$ joining the end-points of $\mu$. \\
Given $D > 0$ there exist $K, \epsilon$ such that if $\mu^h$
penetrates 
exceptional horoballs by at most $D$ then $\mu$ is a $(K, \epsilon )$
quasigeodesic in $\widetilde{M_0}$.
\label{mu-qgeod}
\end{lemma}

\subsection{Ambient and Hyperbolic Geodesics}

\medskip

Recall \\
$\bullet$ $\lambda^h$ $= $ hyperbolic geodesic in $\widetilde{M_{gf}}$
joining $a, b$
\\
$\bullet$ $\lambda = \lambda^a $ $ = $ ambient quasigeodesic in
$\widetilde{M_0}$ constructed from $\lambda^h \subset
\widetilde{M_{gf}}$ \\ 
$\bullet$  $\beta^h$ $ = $ geodesic in $\widetilde{N^h}$ joining $a,
b$ \\
$\bullet$   $\beta_{amb}^0$  $ = $  ambient quasigeodesic 
 in $\widetilde{N_{0}}$ obtained from  $\beta^h$ by replacement of
 hyperbolic by `Euclidean' geodesic segments for $Z$-horoballs in
$\widetilde{N^h}$ \\
$\bullet$ $\beta_{amb} = \Pi_\lambda ( \beta_{amb}^0  )$ \\

\smallskip

By construction, the hyperbolic geodesic $\beta^h$ and the ambient
quasigeodesic  $\beta_{amb}^0$ agree exactly off
horoballs.  $\beta_{amb}$ is constructed from $\beta_{amb}^0$ by
projecting it onto $B_\lambda$ and so by Theorem \ref{mainref},
it is an ambient quasigeodesic. But it might
`backtrack'. Hence, we shall modify it such that it satisfies the no
backtracking condition. First, observe by Theorem \ref{farb} that  $\beta^h$, $\beta_{amb}^0$ track each other off
 some $K$-neighborhood of horoballs.

The advantage of working with $\beta_{amb}$ is that it lies on
$B_\lambda$. 
However, it might backtrack. Now, recall that $B^a( \lambda )$ was
constructed from $B_\lambda$ by replacing each $ B ( \mu )$ by $B^a (
\mu )$. We can therefore choose an ambient quasigeodesic lying on $B^a
( \lambda )$ that tracks $\beta_{amb}$ throughout its length. The
advantage of switching to $B^a ( \lambda )$ is that it is constructed
from ambient geodesics without backtracking lying on the universal
covers of equispaced pleated surfaces. We shall (abusing
notation) call this new ambient quasigeodesic $\beta_{amb}$. Thus
$\beta_{amb}$ now lies on $B^a ( \lambda )$ rather than
$B_\lambda$.

\begin{lemma}
There exists $C > 0$ such that for all
 $x \in \mu^b_i \subset B^a ( \mu , \lambda ) \subset B^a ( \lambda )$
 if $\lambda^h$ lies outside
 $B_n(p)$ for a fixed reference point $p \in \widetilde{M^h}$, then
 $x$ lies outside an $\frac{n-C}{C+1}$ ball about $p$ in
 $\widetilde{N_0}$.
\label{far1}
\end{lemma}

{\bf Proof:} By Corollary \ref{qgray3}, $r_x ( -1 ) \in \lambda^b$.
Since $\lambda^{b}$ is a part of $\lambda^h$, 
therefore $r_x (-1)$ lies outside $B_n (p)$. By Corollary
\ref{qgray3}, there exists $C > 0$ such that for all $i, j \in \{ -1,
0, 1, 2, \cdots$, 

\begin{center}

$|i-j| \leq d(r_x (i) , r_x (j)) \leq C|i-j|$

\end{center}

Also, $d(x, p) \geq i$ since $x \in \mu^b_i$. (Here distances are all
measured in $\widetilde{N_0}$.) Hence,

\begin{eqnarray*}
d(x,p) & \geq & max \{ i, n - C(i+1) \} \\
       & \geq & \frac{n-C}{C + 1}
\end{eqnarray*}

This proves the result. $\Box$

\smallskip

If $x \in B^a (\lambda ) $, then $x \in \lambda^b$ or  $x \in
\lambda^c$ or $x \in B^b ( \mu )$  or $x \in B^c ( \mu )$
for some $\mu $. Hence $x \in B^a ( \lambda )$ implies that either $x$
lies on some horosphere bounding some $ \bf{H} \in \mathcal{H}$ or,
from Lemma \ref{far1} above, $d(x, p) \geq  \frac{n-C}{C + 1}$.
Since $\beta_{amb}$ lies on $B^a ( \lambda )$, we conclude that
$\beta_{amb}$ is an ambient quasigeodesic  in
$\widetilde{N_0}$ such that every point $x$ on $\beta_{amb}$ either
lies on some horosphere bounding some $ \bf{H} \in \mathcal{H}$ or,
from Lemma \ref{far1} above, $d(x, p) \geq  \frac{n-C}{C + 1}$.

McMullen \cite{ctm-locconn} shows (cf Theorem \ref{ctm} )
that in $\widetilde{N^h}$,
any such ambient quasigeodesic
$\beta_{amb}$ lies in a bounded neighborhood of $\beta^h \cup
{ \mathcal{H} } ( \beta^h )$.
We do not as yet know that $\beta_{amb}$ does not backtrack, but we
can convert it into one without much effort. (Note that
Theorem \ref{ctm} does not require "no backtracking".)
Let $\Pi$ denote nearest point projection of  $\widetilde{N^h}$ onto
 $\lambda^h \cup \mathcal{H}$. Since $\Pi$ is a large-scale Lipschitz retract  (Theorem \ref{mainref} ), $\Pi ( \beta_{amb} ) = \beta_1$ is again an
ambient quasigeodesic in $\widetilde{N_0}$. Further, $\beta_1$ tracks
$\beta_{amb}$ throughout its length as $\Pi$ moves points on
$\beta_{amb}$ through
a uniformly bounded distance (Theorem \ref{ctm}). Now $\beta_1$ might backtrack, but it can
do so in a trivial way, i.e. if $\beta_1$ re-enters a horoball after
leaving it, it must do so at exactly the point where it leaves
it. Removing these `trivial backtracks', we obtain an  {\bf
  ambient quasigeodesic without backtracking} $\beta$ which tracks
$\beta_{amb}$ throughout its length. 

\smallskip

{\bf Note:} On the one hand $\beta$ is an ambient quasigeodesic
without backtracking. Hence, it reflects the intersection pattern of
$\beta^h$ with horoballs. On the other hand, it {\bf tracks}
$\beta_{amb}$
whose properties we already know from Corollary \ref{qgray3} above.

\smallskip

Since, of $\beta$ and $\beta^h$, one is an ambient quasigeodesic
without backtracking, and the other a hyperbolic geodesic joining the
same pair of points, we conclude from Theorem \ref{farb} that
they have similar intersection patterns with horoballs, i.e.
there exists $C_0$ such that \\
$\bullet$ If only one of $\beta$ and $\beta^h$ penetrates or travels along the
  boundary of a horoball 
  $\bf{H}$, then it can do so for a distance $ \leq C_0$. \\
$\bullet$ If both $\beta$ and $\beta^h$ enter (or leave) a horoball
  $\bf{H}$ then their entry (or exit) points are at a distance of at
  most $C_0$ from each other.

Again, since $\beta$ tracks $\beta_{amb}$, we conclude that there
 exists
$C > 0$ such that $\beta$ lies in a $C$-neighborhood of $\beta_{amb}$
and hence from Lemma \ref{far1} \\
$\bullet$ Every point $x$ on $\beta$ either
lies on some horosphere bounding some $ \bf{H} \in \mathcal{H}$ or,
$d(x, p) \geq  \frac{n-C}{C + 1} -C$

\smallskip

The above three conditions on $\beta$ and $\beta^h$ allow us to 
deduce the following (identical to the third) condition for $\beta^h$
in the p-incompressible case.

\begin{prop} 
Suppose $(M,P)$ has p-incompressible boundary.
Then every point $x$ on $\beta^h$ either
lies inside some horoball $ \bf{H} \in \mathcal{H}$ or,
$d(x, p) \geq  \frac{n-C}{C + 1} -C = m(n)$
\label{far2}
\end{prop}

We have denoted $ \frac{n-C}{C + 1} -C$ by $ m(n)$, so that $m(n)
\rightarrow \infty$ as $n \rightarrow \infty$.
The above Proposition asserts that the geodesic $\beta^h$ lies outside
large balls about $p$ modulo horoballs. By Lemma \ref{contlemma} this
is almost enough to guarantee the existence of a Cannon-Thurston map.

Again as in Corollary \ref{qgray4} it is not necessary to restrict
ourselves to the p-incompressible case. We deduce, using Corollary
\ref{qgray4} and Lemma \ref{mu-qgeod} the following:

\begin{prop}
Suppose $(M,P)$ has incompressible boundary $\partial_0 P$.
Given $D, n$ there exist $m(n,D)$ such that the
  following holds: \\
Let $\lambda^h$, $\beta^h$, $\lambda$ and $\mu$ be as before. If each $\mu$
penetrates exceptional horoballs by at most $D$,
then every point $x$ on $\beta^h$ either
lies inside some horoball $ \bf{H} \in \mathcal{H}$ or,
$d(x, p) \geq   m(n,D)$, for some function $m(n,D)$
where $m(n,D) \rightarrow \infty$ as $n
\rightarrow \infty$ for each fixed $D$.
\label{far3}
\end{prop}

We can divide $\beta^h$ into two subsets $\beta^c$ and $\beta^b$ as
earlier.
$\beta^c = \beta^h \cap \mathcal{H}$ is the intersection of $\beta^h$
with horoballs, and $\beta^b = \beta^h - \beta^c$.
The main theorem of this paper under the assumption of {\it
  p-incompressibility} follows. Recall that  $M_{gf}$ is
identified with its homeomorphic image in $ N^h $ 
taking cusps to cusps.

\begin{theorem}
Suppose that $N^h \in H(M,P)$ is a hyperbolic structure of bounded geometry
on a pared manifold $(M,P)$ with p-incompressible boundary. Let
$M_{gf}$ denotes a geometrically finite hyperbolic structure adapted
to $(M,P)$, then the map  $i: \widetilde{M_{gf}}
\rightarrow \widetilde{N^h}$ extends continuously to the boundary
$\hat{i}: \widehat{M_{gf}}
\rightarrow \widehat{N^h}$. If $\Lambda$ denotes the limit set of
$\widetilde{M}$, then $\Lambda$ is locally connected.
\label{main}
\end{theorem}

{\it Proof:} Let $\lambda^h$ be a geodesic segment in $\widetilde{M_{gf}}$
lying outside $N_n{(p)}$ for some fixed reference point $p$. Fix 
neighborhoods of the cusps and lift them to the universal cover. Let
$\mathcal{H}$ denote the set of horoballs. Assume without loss of
generality that $p$ lies outside horoballs. Let $\beta^h$ be 
the hyperbolic geodesic in $\widetilde{N^h}$ joining the endpoints of
$\lambda^h$. Further, let $\beta^h = \beta^b \cup \beta^c$ as above.
Then by Proposition \ref{far2}, $\beta^b$ lies outside an $m(n)$ ball
about $p$, with 
 $m(n)
\rightarrow \infty$ as $n \rightarrow \infty$.

Next, let ${\bf{H}}_1$ be some horoball 
that $\beta^h$ meets. Then the entry and exit
points $u$ and $v$ of $\beta^h$ into and out of  ${\bf{H}}_1$
lie  outside an $m(n)$ ball
about $p$. Let $z$ be the point on the boundary sphere that
${\bf{H}}_1$ is based at. Then for any sequence  $x_i \in
{\bf{H}}_1$
with $d(p, x_i ) \rightarrow \infty$, $x_i \rightarrow z$. 
If $\{ x_i \}$ and $\{ y_i \}$ denote two such sequences, then the
visual diameter of the set $\{ x_i , y_i \}$ must go to zero. Hence,
if $[x_i,y_i]$ denotes the geodesic joining $x_i , y_i$ then
$d(p,[x_i,y_i]) \rightarrow \infty$. Since, $u, v$ lie outside an
$m(n)$ ball, there exists some function $\psi$, such that the geodesic
$[u,v]$
lies outside a $\psi (m(n))$ ball around $p$, where $\psi (k)
\rightarrow \infty$ as $k \rightarrow \infty$. 

We still need to argue that the function $\psi$ is independent
of the horoball ${\bf{H}}_1$. Assume that
the lifts of horoballs miss the base-point $p$. If the conclusion fails,
there exist a sequence of horoballs ${\bf{H}}_i$, points $x_i, y_i$
on the boundary of ${\bf{H}}_i$, such that $[x_i, y_i]
\subset {\bf{H}}_i$ passes
through some fixed $N$-ball $B_N(p)$ about $p$. Passing to a subsequence
and extracting a limit we may assume that the
horoballs ${\bf{H}}_i$ converge to some horoball ${\bf{H}}$
cutting  $B_N(p)$
(possibly touching $p$). Then the above discussion for 
${\bf{H}}_1$ furnishes the required contradiction.
Hence, the choice of this function $\psi$ does not depend on
 ${\bf{H}}_1$. We conclude that there
exists such a function for all of $\beta^c$. We have thus established:
 \\
$\bullet$ $\beta^b$ lies outside an $m(n)$ ball about $p$. \\
$\bullet$ $\beta^c$ lies outside a $\psi (m(n))$ ball about $p$. \\
$\bullet$ $ m(n)$ and  $\psi (m(n))$ tend to infinity as $n \rightarrow
 \infty$\\

Define $f(n) = min( m(n),\psi (m(n)) )$. 
Then  $\beta^h$ lies outside an $f(n)$ ball about $p$ and 
 $f(n)
\rightarrow \infty$ as $n \rightarrow \infty$.

 By Lemma \ref{contlemma}  $i: \widetilde{M_{gf}}
\rightarrow \widetilde{N^h}$ extends continuously to the boundary
$\hat{i}: \widehat{M_{gf}}
\rightarrow \widehat{N^h}$. This proves the first statement of the
theorem.

Now, for any geometrically finite Kleinian group, its limit set is
locally connected.( See, for instance, \cite{and-mask}.) Hence, the
limit set of $\widetilde{M_{gf}}$ is locally connected.  Further, the
continuous image of a compact locally connected set is locally
connected \cite{hock-young}. Hence, if $\Lambda$ denotes the limit set of
$\widetilde{N^h}$, then $\Lambda$ is locally connected. This proves
the theorem. $\Box$

Again as in Corollary \ref{qgray4} and Proposition \ref{far3}
it is not necessary to restrict
ourselves to the p-incompressible case. In Theorem \ref{main} above we
removed the clause on being within horoballs from
Proposition \ref{far2}. In the proof of Theorem \ref{main} above the
one new thing we had shown was that geodesics entering and leaving
horoballs
at a large distance from the reference point $p$, itself lies at a
large distance from $p$. Using this observation along with 
Proposition \ref{far3} we easily  deduce:

\begin{cor}
Suppose $(M,P)$ has incompressible boundary $\partial_0 P$.
Given $D, n$ there exist $m(n,D)$ such that the
  following holds: \\
Let $\lambda^h$, $\beta^h$  be as before. If  $\lambda^h$
penetrates exceptional horoballs by at most $D$
and if $\lambda^h$ lies outside $B_n (p)$ in $\widetilde{M_{gf}}$
then  $\beta^h$ lies outside a ball of radius
$  m(n,D)$, for some function $m(n,D)$
where $m(n,D) \rightarrow \infty$ as $n
\rightarrow \infty$ for each fixed $D$.
\label{maincor}
\end{cor}

The above proposition will be useful in the next subsection when we go
from
{\it p-incompressibility} to {\it incompressibility}.

\subsection{From p-incompressibility to Incompressibility}

In this subsection we shall use Corollary \ref{maincor} to relax the
hypothesis of p-incompressibility. Recall that we started with a
hyperbolic geodesic $\lambda^h$ in $\widetilde{M_{gf}}$ and then
modified it along horoballs to obtain $\lambda$. In this subsection we
first consider subsegments $\lambda^h_1  , \cdots \lambda^h_k$ where
each $\lambda^h_i$ starts and ends on horospheres bounding horoballs,
and the complementary segments (of $\lambda^h$) lie inside
horoballs. This decomposition is made in such a way that the starting
and ending points of $\lambda^h_i$ lie on {\it exceptional horoballs}
(See previous Subsection for definitions.) Having completed this
decomposition we consider hyperbolic geodesics $\beta^h_i$ in
$\widetilde{N^h}$ joining the
end-points of $\lambda^h_i$. If we ensure that the geodesics
$\lambda^h_i$ penetrate  exceptional horoballs 
by a uniformly bounded amount, then by Corollary \ref{maincor} above, we can
ensure that the segments $\beta^h_i$ lie outside a large ball about
$p$.

If further, we can ensure that \\
$\bullet$ the $\beta^h_i$ penetrate the horoballs
they start and end on by a bounded amount \\
$\bullet$ and that the terminal point
of $\beta^h_i$ and the initial point of $\beta^h_{i+1}$ are separated
by more than a critical amount,\\
 then (by appealing to the fact that
local quasigeodesics are global quasigeodesics \cite{gromov-hypgps}
\cite{GhH}) we conclude that the union of the $\beta^h_i$ with the
hyperbolic geodesics interpolating between them (and lying entirely
within horoballs) is a (uniform) hyperbolic quasigeodesic.

Since each $\beta^h_i$ lies outside a large ball about $p$, it follows
that their union along with interpolating geodesics does so too.
This union being a quasigeodesic, the geodesic $\beta^h$ joining the
end-points of $\lambda$ must also lie outside a large ball about
$p$. By Lemma \ref{contlemma} the existence of a Cannon-Thurston map
follows.

We now furnish the details of the above argument.

\smallskip

First note from Corollary \ref{maincor} that given $D$ there exists
$m(n,D)$
such that if $\lambda^h_1$ be a
subsegment of $\lambda^h$ penetrating exceptional horoballs by at most
$D$, and if $\lambda^h$ (and hence $\lambda^h_1$) lies outside
$B_n(p)$
then $\beta^h_1$ (the hyperbolic geodesic in $\widetilde{N^h}$ joining
the endpoints of $\lambda^h_1$) lies outside an $m(n,D)$-ball about $p$.

The next proposition follows from the fact that local quasigeodesics
are global quasigeodesics in hyperbolic space. (See Gromov
\cite{gromov-hypgps} p.187, Prop. 7.2C, \cite{CDP}.)

\begin{prop}
There exist $D, K, \epsilon $ such that the following holds:
\\
Suppose $\beta^x$ is a path in $\widetilde{N^h}$ such that $\beta^x$
can be decomposed into finitely many  geodesic segments $\beta_1,
\cdots \beta_k$. Further suppose that the starting and ending point of
each $\beta_i$ lie on  exceptional horospheres (except possibly the
starting point of $\beta_1$ and the ending-point of $\beta_k$)
meeting the horospheres at right angles. Also
suppose that the `even segments' $\beta_{2i}$ lie entirely within
exceptional horoballs and have length greater than $D$. Then $\beta^x$
is a $(K, \epsilon )$ quasigeodesic.
\label{local-global}
\end{prop}

Proof Idea: Since geodesic segments lying outside horoballs
meet horospheres at right angles, 
successive pieces meet at an angle that is bounded below.
Hence, if $x_i, z_i \in \beta_i , \beta{i+1}$ respectively
and $\beta_i \cap \beta{i+1} = y_i$, then $[x_i,y_i] \cup
[y_i,z_i]$ is a (uniform) $(K_1, \epsilon_1 )$-quasigeodesic
of length $\geq D+1$ assuming that any two horoballs are
separated by a distance of at least one. That is every
segment of length $ D+1$ is a $(K_1, \epsilon_1 )$-quasigeodesic,
i.e. (See \cite{CDP} Ch. 3) it is a local quasigeodesic.
If $D$ is sufficiently long, stability of quasigeodesics
(Theorems 1.2, 1.3, 1.4 of  \cite{CDP} Ch. 3)
ensures that $\beta^x$
is a $(K, \epsilon )$ quasigeodesic for some $K, \epsilon$ independent
of $\beta$.

\smallskip

In fact we do not need that the horoballs be exceptional, but that
they be uniformly separated, i.e.  the distance between any two
horospheres is uniformly bounded below. This is clearly true for the
exceptional horoballs since they are lifts of certain given cusps
in $M_{gf}$ which in turn are finitely many in number and have been
chosen  uniformly separated.

Next, we need to look closely at how we use Corollary \ref{qgray3} to
prove Theorem \ref{main}. Note that $\beta_{amb}$ is an ambient
quasigeodesic independent of the hypothesis of {\it
  p-incompressibility}. It is in concluding that $\beta_{amb}$ lies
outside a large $m(n)$-ball about the reference point $p$ that we
needed to construct quasigeodesic rays going down to $\lambda$.

Now let $\lambda^h = [a,x]$ be a hyperbolic geodesic in
$\widetilde{M_{gf}}$ having the end-point $x$ on an exceptional
horosphere bounding a horoball {\bf $H$}. By our construction of
$B_\lambda$ we note that $B_\lambda$ is quasi-isometrically embedded
in $\widetilde{N}$ and further that $\beta_{amb}$ (the ambient
quasigeodesic corresponding to $\lambda^h$ ) meets {\bf $H$} either at
$x$ (if {\bf $H$} is a $(Z+Z)$-horoball) or at some point on 
$r_x$ (for a quasigeodesic ray $r_x$ constructed through $x$ and lying
in the universal cover $\widetilde{E}$ of an end $E$ of the manifold) 
if {\bf $H$} be a $Z$-horoball.

Now suppose $\lambda^h_1 = [a,x]$, $\lambda^h_2 = [x,y]$ and
$\lambda^h_3 = [y,b]$ be three hyperbolic geodesic segments such that
$[x,y] \subset {\bf H}$ and $\lambda^h_3$ meets {\bf $H$} at
$y$. Again, as before the corresponding ambient quasigeodesic
$\beta_{amb}^1$ meets ${\bf H}$ at $x$ or at some point along an $r_x$
and  $\beta_{amb}^3$ meets ${\bf H}$ at $y$ or at some point along an
$r_y$. 

In any case, the distance between the entry point of
$\beta^h_1$ (the hyperbolic geodesic joining the end-points of
$\beta^1_{amb}$) into {\bf $H$} and the exit point of $\beta^h_3$ (the
hyperbolic geodesic joining the end-points of $\beta^3 _{amb}$)
from {\bf $H$} is greater than $d(x,y) - 2C_0$ for some $C_0$
depending on the quasiconvexity constant of $\beta^i_{amb}$ (i=1,3) by
Theorem \ref{farb}. Thus $C_0$ depends only on the quasi-Lipschitz
constant of $\Pi_\lambda$. 

\smallskip

We are now in a position to break $\lambda^h$ into pieces. Let
$\lambda^h_{2i}$ denote the maximal subsegments of $\lambda^h$ lying
inside exceptional horoballs and having length greater than $(2C_0 +
D_0)$. Let $\lambda^h_{2i-1}$ denote the complementary segments. Now,
let $\overline{\beta^h_{2i-1}}$ be the hyperbolic geodesic in
$\widetilde{N^h}$ joining the endpoints of  ${\lambda^h_{2i-1}}$. Then
the entry point of  $\overline{\beta^h_{2i-1}}$ into the {\it
  exceptional horoball} {\bf H}  it terminates on lies at a distance
greater than $(2C_0 + D_0 - 2C_0)=D_0$ from the exit point of
$\overline{\beta^h_{2i+1}}$ from the same horoball {\bf H}. 
Shorten the hyperbolic geodesic  $\overline{\beta^h_{2i-1}}$ if
necessary by cutting it off at the entry point into {\bf H}. Let 
${\beta^h_{2i-1}}$ denote the resultant geodesic. By Corollary
\ref{maincor}
 ${\beta^h_{2i-1}}$ being a subsegment of  $\overline{\beta^h_{2i-1}}$
lies at a distance of at least $m(n,(2C_0 + D_0)) = m_1(n)$ from the reference
point $p$.

Next denote by  ${\beta^h_{2i}}$ the hyperbolic geodesic lying
entirely within {\bf H} joining the entry point of
${\beta^h_{2i-1}}$ into {\bf H} to the exit point of
${\beta^h_{2i+1}}$ from {\bf H}. The initial and terminal points of
$\beta^h_{2i}$ lying on {\bf H} are at a distance of at least $m_1(n)$
from $p$. Therefore each $\beta^h_{2i}$ and hence the union of all the
$\beta^h_i$ lie outside an $m_2(n)$ ball about $p$ where $m_2(n)
\rightarrow \infty$ as $n \rightarrow \infty$. Further, the union of
the segments $\beta^h_i$ is a hyperbolic quasigeodesic by Proposition
\ref{local-global} and hence lies at a bounded distance $D^{\prime}$
from the hyperbolic geodesic $\beta^h$ joining the endpoints of
$\lambda$ (i.e. the end-points of $\lambda^h$). Let $m_3(n) = m_2(n) -
D^{\prime}$. Then what we have shown amounts to : \\
$\bullet$ If $\lambda^h$ lies outside $B_n(p)$, then $\beta^h$ lies
outside a ball of radius $m_3(n)$ about $p$. \\
$\bullet$  $m_3(n)
\rightarrow \infty$ as $n \rightarrow \infty$ \\

Coupled with Lemma \ref{contlemma} this proves the  main
theorem of this paper given below:

\begin{theorem}
Suppose that $N^h \in H(M,P)$ is a hyperbolic structure of bounded geometry
on a pared manifold $(M,P)$ with incompressible boundary $\partial_0 M$. Let
$M_{gf}$ denotes a geometrically finite hyperbolic structure adapted
to $(M,P)$. Then the map  $i: \widetilde{M_{gf}}
\rightarrow \widetilde{N^h}$ extends continuously to the boundary
$\hat{i}: \widehat{M_{gf}}
\rightarrow \widehat{N^h}$. If $\Lambda$ denotes the limit set of
$\widetilde{M}$, then $\Lambda$ is locally connected.
\label{main2}
\end{theorem}

\section{Examples and Consequences}

As mentioned in the Introduction, the simplest non-trivial examples to
which Theorem \ref{main2} applies are hyperbolic three manifolds of
finite volume fibering over the circle. These include the original
examples of Cannon and Thurston \cite{CT} as well as those of Bowditch
\cite{bowditch-ct} \cite{bowditch-stacks}. 

The next set of examples are those three manifolds homeomorphic to the
product of a surface and $\Bbb{R}$. These were dealt with in Minsky's
work \cite{minsky-jams} and the punctured surface case was dealt with by
Bowditch \cite{bowditch-ct} \cite{bowditch-stacks}.

The case of three manifolds of freely indecomposable fundamental group
were dealt with independently by Klarreich\cite{klarreich} and  the
author \cite{mitra-trees}. In fact, Klarreich's theorem is really a
reduction theorem which effectively says that if one can prove
Cannon-Thurston for closed surface groups (of some given geometry)
then one can also prove it for 3 manifolds whose ends have the same
geometry. In combination with Minsky's adaptation of the original
proof of Cannon-Thurston, this proves the theorem for bounded geometry
3-manifolds with incompressible (closed surface) boundary. 
In \cite{mitra-trees} we had approached this problem  directly
 and had given different proofs of these results. In this paper we have
continued the approach in \cite{mitra-trees} to prove the analogue of
the above results in the presence of parabolics. 

{\bf Incompressibility and compressibility}

One problem that has not been addressed by this paper or its
predecessor \cite{mitra-trees} is the case of compressible $\partial_0
M$. Miyachi \cite{miyachi-ct} (see also Souto \cite{souto-ct} )
has solved this problem in the bounded
geometry case, when there are no cusps. In \cite{mahan-split} we shall combine the reduction techniques of this paper with a coarse version
of Miyachi's argument to settle affirmatively the question of 
existence of Cannon-Thurston maps for finitely generated
Kleinian groups. 

{\bf Coarse Framework}\\
 Many of the arguments of this paper rightfully belong to
the domain of `coarse' or `asymptotic' or `large-scale' geometry in
spirit and it is more than likely that they may be generalized to this
setting. In \cite{mitra-ct} and \cite{mitra-trees}, we had proven the
following theorems:

\begin{theorem}
\cite{mitra-ct} Let $G$ be a hyperbolic group and let $H$ be a hyperbolic subgroup
that is normal in $G$. Let 
$i : \Gamma_H\rightarrow\Gamma_G$ be the continuous proper
 embedding of $\Gamma_H$ in $\Gamma_G$ described above.
 Then $i$ extends to a continuous
map $\hat{i}$ from
$\widehat{\Gamma_H}$ to $\widehat{\Gamma_G}$.
\label{mitra-ct}
\end{theorem}

\begin{theorem}
\cite{mitra-trees}  Let (X,d) be a tree (T) of hyperbolic
 metric spaces satisfying the
qi-embedded condition.  Let $v$ be a vertex of $T$. If $X$ is hyperbolic
${i_v} : X_v \rightarrow X$ extends 
continuously to ${\hat{{i_v}}} : \widehat{X_v}
\rightarrow \widehat{X}$.
\label{mitra-trees}
\end{theorem}

In this paper we have described a fairly general way of handling
cusps. The generalization to `coarse geometry' involves dealing with
relatively hyperbolic groups {\it a la} Gromov \cite{gromov-hypgps},
Farb \cite{farb-relhyp} and Bowditch \cite{bowditch-relhyp}. 

In \cite{mahan-pal}, the author and Abhijit Pal generalize Theorem
\ref{mitra-trees} to the context of trees of (strongly) relatively hyperbolic
metric spaces. In \cite{pal-ct}, Pal adapts the
techniques of this paper and generalizes Theorem \ref{mitra-ct}
to (strongly) relatively hyperbolic  normal subgroups
of (strongly) relatively hyperbolic  groups.

\bibliography{pared}
\bibliographystyle{alpha}

\end{document}